\documentclass[11pt,a4paper]{amsart}

\usepackage{amssymb}
\usepackage{amsmath}
\usepackage{amsfonts}
\usepackage{amsthm}
\usepackage{verbatim}

\renewcommand{\phi}{\varphi}
\renewcommand{\epsilon}{\varepsilon}
\renewcommand{\hat}{\widehat}
\renewcommand{\tilde}{\widetilde}

\newcommand{\C}{\mathbb{C}}

\newcommand{\HH}{\mathcal{H}}
\newcommand{\PP}{\mathcal{P}}
\newcommand{\VV}{\mathcal{V}}

\newcommand{\Ub}{\overline{U}}
\newcommand{\Xb}{\overline{X}}
\newcommand{\Bbar}{\overline{B}}
\newcommand{\Wh}{\hat{W}}

\newcommand{\hhat}{\tilde{h}}
\newcommand{\hhatinv}{(\tilde{h}^{-1})}
\newcommand{\ginv}{(g^{-1})}
\newcommand{\xit}{\tilde{\xi}}
\newcommand{\that}{\hat{\theta}}
\newcommand{\ehat}{\hat{\epsilon}}
\newcommand{\phat}{\hat{\phi}}
\newcommand{\phicy}{\rho}
\newcommand{\xitbar}{\overline{\xit}}

\newcommand{\hook}{\lrcorner\:}
\newcommand{\zbar}{\overline{z}}
\newcommand{\wbar}{\overline{w}}

\newcommand{\ddbar}{\partial\overline{\partial}}
\newcommand{\dbar}{\overline{\partial}}

\newcommand{\lap}{\Delta_{g_+}}

\newcommand{\varthetat}{\tilde{\vartheta}}
\newcommand{\varthetah}{\hat{\vartheta}}

\newcommand{\phih}{\hat{\phi}}

\newcommand{\al}{\alpha}
\newcommand{\bbar}{\overline{\beta}}
\newcommand{\abar}{\overline{\alpha}}
\newcommand{\gbar}{\overline{\gamma}}
\newcommand{\mbar}{\overline{\mu}}

\newcommand{\delbar}{\overline{\delta}}
\newcommand{\oneb}{\overline{1}}
\newcommand{\zerob}{\overline{0}}

\newcommand{\jb}{\overline{j}}
\newcommand{\kb}{\overline{k}}
\newcommand{\lb}{\overline{l}}
\newcommand{\mb}{\overline{m}}

\newcommand{\Vol}{\textnormal{Vol}}
\newcommand{\const}{\textnormal{const.\:}}
\renewcommand{\Re}{\textnormal{Re}}
\renewcommand{\Im}{\textnormal{Im}}
\newcommand{\Ric}{\textnormal{Ric}}
\newcommand{\Scal}{\textnormal{Scal}}

\newcommand{\Ab}{\boldsymbol{A}}
\newcommand{\Bb}{\boldsymbol{B}}
\newcommand{\Db}{\boldsymbol{\Delta}}
\newcommand{\Deeb}{\boldsymbol{D}}
\newcommand{\Mb}{{\boldsymbol M}}
\newcommand{\Ubo}{\boldsymbol{U}}
\newcommand{\Xbo}{\boldsymbol{X}}
\newcommand{\nb}{n^\prime}
\newcommand{\gb}{\boldsymbol{g}}
\newcommand{\gbp}{\boldsymbol{g_+}}
\newcommand{\gbt}{\boldsymbol{\tilde{g}}}
\newcommand{\pb}{\psi}

\theoremstyle{plain}
\newtheorem{thm}{Theorem}[section]
\newtheorem{prop}[thm]{Proposition}
\newtheorem{lem}[thm]{Lemma}
\newtheorem{cor}[thm]{Corollary}

\theoremstyle{definition}

\theoremstyle{remark}
\newtheorem{rem}{Remark\!\!}

\theoremstyle{plain}

\theoremstyle{plain}

\begin{document}
 
\title[Volume renormalisation]{Volume renormalization for complete 
Einstein--K\"ahler metrics}

\author{Neil Seshadri}

\address{Graduate School of Mathematical Sciences, The University of Tokyo
        \newline\indent 3--8--1 Komaba, Meguro, Tokyo 153--8914, Japan}
\email{seshadri@ms.u-tokyo.ac.jp}

\begin{abstract}
For a strictly pseudoconvex domain in a complex manifold we define a 
renormalized volume with respect to the approximately Einstein complete 
K\"ahler metric of Fefferman. We compute the conformal anomaly in complex 
dimension two and apply the result to derive a renormalized 
Chern--Gauss--Bonnet formula.
Relations between renormalized volume and the CR $Q$-curvature are also 
investigated.
\end{abstract}

\maketitle

\numberwithin{equation}{section}

\section{Introduction}
There has been much recent activity in the area of 
conformal geometry centred around volume renormalisation. In this note we 
demonstrate that a similar notion of 
volume renormalisation exists in CR geometry and obtain some explicit 
formulae, including a 
renormalised Chern--Gauss--Bonnet formula, in complex dimension two.

The quantity in conformal geometry known as renormalised volume was first 
defined by 
Henningson--Skenderis~\cite{hs} and Graham~\cite{graham} in their studies
of conformally compact Einstein manifolds. A conformally compact (Einstein) 
manifold consists of a compact manifold-with-boundary $\Xb$ and an (Einstein) 
metric $g_+$ on the interior $X$ satisfying the condition that 
$\phi^2g_+$ extends 
at least continuously to the boundary $M$, where $\phi$ is a defining function
 for $M$, i.e., $\phi < 0$ on $X$ and $\phi = 0, d\phi \neq 0$ on $M$. The 
restriction of $\phi^2g_+$ to $TM$ rescales upon changing 
$\phi$ whence defining a conformal class of metrics $[g]$ on $M$. 
Fefferman--Graham~\cite{fg85} called $g_+$ the 
Poincar\'e metric associated to the conformal manifold $(M,[g])$. The volume 
of 
$(X,g_+)$ is infinite; but one can renormalise it by considering the expansion,
in powers of a special defining function, of the volume of the set 
$\{\phi < \epsilon\}$ as $\epsilon\to 0$. The 
constant term in this asymptotic expansion, called the renormalised volume, 
is an invariant of the metric $g_+$, provided $M$ is odd-dimensional. 
If $M$ is even-dimensional, the renormalised volume is dependent on the 
choice of representative $g$; however in this case there is a log term in the
volume expansion whose coefficient turns out to be an invariant of $g_+$. The
 difference between the renormalised volumes corresponding to two different 
choices of conformal representative in the even-dimensional case is known as
 the conformal anomaly. It is the boundary integral of a nonlinear partial
 differential operator applied to the conformal scale.

To explore an analogous notion in 
CR geometry, let $M$ be a compact 
strictly pseudoconvex $(2n+1)$-dimensional 
CR manifold embedded into a complex $(n+1)$-manifold $X$, with 
natural CR structure coming from the embedding. 
There is conformal class of pseudohermitian Levi metrics on $M$, the 
``pseudo'' coming from the fact that they only act on a certain subbundle of
$TM$. Each Levi metric in this conformal class corresponds to a contact 
one-form called a pseudohermitian structure. Consider the complete K\"ahler 
metric $g_+$ with K\"ahler form
\begin{equation*}
\omega = \partial\overline{\partial}\log\left(-\frac{1}{\phicy}\right)
       = -\frac{\ddbar\phicy}{\phicy} 
+ \frac{\partial \phicy\wedge\dbar \phicy}{\phicy^2},
\end{equation*}
for any sufficiently smooth defining function $\rho$.
If $\HH_\rho$ denotes the bundle whose fibre over each 
$M_\rho^\epsilon := \{\phicy = \epsilon\}$ is $\HH_\rho^\epsilon$,
the holomorphic tangent bundle to $M_\rho^\epsilon$, then
$\phicy g_+ |_{\HH_\rho^\epsilon}$ extends continuously to $M$ and 
$$
\phicy g_+ |_{\HH_\rho^0} \in (\textnormal{conformal class of Levi metrics on
 $M$}).
$$
Thus $g_+$ provides a CR-geometric analogue of a conformally compact 
metric. Fefferman~\cite{fef76} produced a $\rho$ that locally solves 
a certain Monge--Amp\'ere type initial-value problem. For this $\rho$, $g_+$ 
is approximately Einstein near $M$, i.e., satisfies
$$
\textnormal{Ricci form of $g_+$} = -(n+2)\omega + O(\rho^{n+1}).
$$
Remark that in the
special case where $X$ is a Stein manifold, the existence and uniqueness of a
complete Einstein--K\"ahler metric, asymptotic to Fefferman's approximately
Einstein metric, is given by Cheng--Yau~\cite{cy}.
By applying techniques from the conformal 
setting we obtain an asymptotic expansion for Fefferman's approximately 
Einstein metric in powers of
 a special defining function $\phi$ with coefficients pseudohermitian 
invariants of the chosen pseudohermitian representative.  We can
extend Fefferman's approximately Einstein metric to an Hermitian metric over
the whole of $X$ and then consider the asymptotic expansion of the volume of 
the set $\{\phi < \epsilon\}$: 
\begin{equation}
\label{eq:volume}
c_0\epsilon^{-n-1} + c_1\epsilon^{-n}
+\cdots + c_n\epsilon^{-1} + L\log(-\epsilon) + V + o(1).
\end{equation}
The constant term $V$ is defined to be the renormalised volume. The 
coefficients $c_j$ and $L$ are integrals over $M$ of local pseudohermitian 
invariants of the chosen pseudohermitian representative, with $L$ being 
independent of this choice. While we shall illustrate the volume 
renormalisation procedure in general dimension, it turns out that obtaining 
explicit formulae seems to be a computationally infeasible task (by hand) 
except in the lowest dimension. Restricting attention to the $n=1$ case 
gives our main result.

\begin{thm} 
\label{thm:main}
Let $M$ be a smooth compact strictly pseudoconvex 
CR manifold of dimension three that forms the boundary of a complex manifold 
$X$. Let $g_+$ be Fefferman's approximately Einstein metric on $X$ 
near $M$. Let $\theta$ and $\that = e^{2\Upsilon}\theta$ be two 
pseudohermitian structures on $M$. Then the 
conformal anomaly of the renormalised volumes $V_\theta, V_{\that}$ of 
$(X,g_+)$ with respect to $\theta, \that$, respectively, is
\begin{equation*}
\begin{split}
V_{\that} - V_\theta &= \int_M \frac{1}{96}\left(\Delta_b\Scal 
- 8\Im A_{11,}^{\phantom{11,}11} \right)\Upsilon 
+ \frac{1}{16}\Scal\Upsilon_1\Upsilon^1\\
&\quad - \frac{1}{4}\Big((\Upsilon_1\Upsilon^1)^2 
+ \Im(\Upsilon_1\Upsilon_1A^{11})
+ \frac{1}{4}(\Upsilon_T)^2\\ 
&\quad - 2\Re(\Upsilon_{11}\Upsilon^1\Upsilon^1) 
+ \Upsilon_1\Upsilon^1\Delta_b\Upsilon  \Big) \: \theta\wedge d\theta .
\end{split}
\end{equation*}
Here $\Scal$ and $A_{11}$ denote the pseudohermitian scalar curvature and 
torsion of $(M,\theta)$, $\Delta_b$ is the sublaplacian and the subscripts $1$
and $T$ denote pseudohermitian covariant differentiation in the directions of 
the CR structure and characteristic vector field, respectively. Indices are 
raised with the Levi metric of $\theta$.

Furthermore, in the volume expansion (\ref{eq:volume}), $L=0$.
\end{thm}

The proof of Theorem~\ref{thm:main} rests largely on computing the asymptotic 
expansion for $g_+$ in the form mentioned above. The tool for this 
purpose is the ambient connection of Graham--Lee~\cite{gl}, which is an 
extension of the pseudohermitian connection inside the boundary. The natural
coframe used in~\cite{gl} is one that restricts to an admissible coframe on 
each pseudohermitian manifold of a foliation near the boundary. However, with 
the strategy of simplifying Einstein equation calculations by putting $g_+$ 
into a diagonal form, we find that we must modify this frame; that constitutes
the main technical step.  

As a rather straightforward consequence, via some basic invariant theory, of 
Theorem~\ref{thm:main} we prove a renormalised Chern--Gauss--Bonnet formula.

\begin{cor} 
\label{cor:main}
Let $M$ be a smooth compact strictly pseudoconvex 
CR manifold of dimension three that forms the boundary of a complete 
Einstein--K\"ahler manifold $(X,g_{\textnormal{\tiny{EK}}})$. Then 
$$
\chi(X) = \int_X \left(c_2 - \frac{1}{3}c_1^{\phantom{1}2}\right) + \VV(M),
$$
where 
$$
\VV(M) := \frac{1}{\pi^2}\left( 6V - \frac{1}{128}
\int_M (\Scal)^2 - 16|A|^2\ \theta \wedge d\theta\right)
$$
is a global CR invariant.
\end{cor}

Anderson~\cite{an} has obtained an analogous formula for conformally 
compact four-manifolds. Of course in that setting there is no
boundary error term since the renormalised volume itself is conformally 
invariant. Now if the CR manifold $M$ is the boundary of a domain in 
$\C^{n+1}$ then
 from work of Burns--Epstein~\cite{be} we have that 
$$
\chi(X) = \int_X \left(c_2 - \frac{1}{3}c_1^{\phantom{1}2}\right) - \mu(M),
$$
where $\mu(M)$ is the Burns--Epstein global CR invariant of~\cite{be88}. So 
in this case, $\VV(M) = -\mu(M)$.

We remark that Herzlich~\cite{herzlich} (see also~\cite{bh}) has recently 
defined a renormalised volume and obtained a renormalised 
Chern--Gauss--Bonnet formula 
for real four-dimensional generalisations of 
our approximately Einstein complete K\"ahler manifolds, called almost complex 
hyperbolic Einstein (ACHE) manifolds. His formulation, though, differs from 
ours somewhat\footnote{The main difference being in the choice of defining 
function. He also uses a different constant in the 
Einstein equation.} and as a 
consequence the 
explicit formulae obtained there do not coincide with ours. 

The special defining function we use, namely, the unique solution to the 
boundary-value problem
$$
\frac{i}{2}(\dbar\phi - \partial\phi)|_M = \theta\ ;\quad
|\partial(\log(-\phi))|_{g_+}^2 = 1,
$$
has its motivation in the relation provided by Fefferman~\cite{fef76} between 
a CR manifold $M$ and an even-dimensional conformal structure living on a 
canonical circle bundle over $M$. In fact our defining function $\phi$ is the
(negative of the) square of the special defining function associated to 
Fefferman's conformal manifold. 

Fefferman's relation between CR geometry and conformal geometry also provided 
the means (\cite{fh}) for defining the CR $Q$-curvature 
$Q^{\textnormal{CR}}_\theta$ from that concept
in conformal geometry. In general even-dimensional conformal geometry, 
Fefferman--Graham~\cite{fg} showed that the integral of the $Q$-curvature of 
the boundary is a constant multiple of the 
coefficient of the log term in the volume expansion of the interior manifold
with respect to the Poincar\'e metric. Analogously, in the CR case we shall see
that
\begin{equation}
\label{eq:CRQ}
\int_M Q^{\textnormal{CR}}_\theta = \const L.
\end{equation}
It was shown by Fefferman--Hirachi~\cite{fh} that when $n=1$, 
$\int_M Q^{\textnormal{CR}}_\theta
= 0$. So (\ref{eq:CRQ}) verifies the last statement in Theorem~\ref{thm:main}.
It is worth mentioning, however, that although the statment regarding $L$ in 
Theorem~\ref{thm:main} is proved by integrating a divergence term 
(see \S~\ref{sec:n=1}), this divergence term is \textit{not} a constant 
multiple of $Q^{\textnormal{CR}}_\theta$, since it does not 
satisfy the appropriate 
(see~\cite{fh}) transformation law under change in pseudohermitian structure.

The analogy between  
 CR geometry and conformal geometry breaks down when one considers the 
obstruction to determining the high-order expansions for the respective 
metrics $g_+$. It was shown in~\cite{fef76} that the obstruction to 
obtaining a smooth complete Einstein--K\"ahler metric is a scalar function. We 
verify this fact by showing that the trace, with respect to the chosen 
representative Levi metric, of the $\phi^{n+1}$ coefficient in 
the asymptotic expansion for $g_+$ is formally \textit{un}determined, while 
the trace-free part of this coefficient is formally determined. 
This is in contrast to the even-dimensional conformal setting 
(\cite{graham}) where the obstruction first
occurs in determining the trace-free part of the conformally compact metric.

\textit{Layout of the paper.} In \S~\ref{sec:pseudo} we review the basics of 
pseudohermitian geometry and introduce the ambient connection of Graham--Lee. 
In \S~\ref{sec:einstein} we fix a special defining function and suitable
coframe, compute the Levi-Civita connection and Ricci form and describe,
in general dimension, how to determine our desired expansion for $g_+$ from 
the Einstein equation. In \S~\ref{sec:volume} we
define the renormalised volume and conformal anomaly and show that the log
term coefficient is an invariant of Fefferman's approximately Einstein 
metric. The proof of Theorem~\ref{thm:main} is the content of 
\S~\ref{sec:n=1} and Corollary~\ref{cor:main} is proved in \S~\ref{sec:cgb}.
It is illustrative to see how the whole procedure works for a concrete 
example, so we devote \S~\ref{sec:bergman} to explicit calculations for the 
Bergman metric on the unit ball in $\C^2$.
We confine to the Appendix a discussion of CR $Q$-curvature and the proof of
formula (\ref{eq:CRQ}).

\textit{Notations.} Lowercase Greek indices run from 1 to $n$ and, unless 
otherwise indicated, lowercase Latin indices run from 0 to $n$. The letter 
$i$ will denote the quantity $\sqrt{-1}$. We observe the summation convention.

\textit{Acknowledgements.}
This work was part of the author's Master's thesis at the University 
of Tokyo. I am deeply indebted to my supervisor Professor Kengo Hirachi for  
suggesting this particular problem to work on and for his constant support 
and advice. I also thank Professor Robin Graham for useful conversations and 
the referee for suggestions for improvement to an earlier draft. Part of this 
project was carried out at the University of Adelaide, whose mathematics 
department I thank for their hospitality. 
Financial support through the Japanese Ministry of Education 
research scholarship is gratefully acknowledged.

\section{Pseudohermitian geometry}
\label{sec:pseudo}

\subsection{Pseudohermitian manifolds} Recall some standard notions from
pseudohermitian geometry, as contained in,  e.g.,~\cite{lee88}. A 
\textit{CR structure} 
(of hypersurface type) on a real $(2n+1)$-dimensional manifold $M$ is an 
$n$-dimensional complex subbundle, also called the \textit{holomorphic tangent 
bundle}, $\HH^M\subset\C TM$ 
satisfying $\HH^M\cap\overline{\HH^M}=\{0\}$ and
$[\HH^M,\HH^M]\subset\HH^M$.
  A global nonvanishing section $\theta$ of the 
real line bundle $E\subset T^{\ast}M$ of annihilators of 
$H^M := \textnormal{Re }\HH^M\oplus\overline{\HH^M}$ exists
if and only if $E$ is orientable; since $H^M$ is 
naturally oriented by its complex structure, such a $\theta$ exists if and 
only if $M$ is orientable. If a $\theta$ exists then clearly it is 
determined by the CR structure only up to conformal multiple.
If the Hermitian bilinear 
\textit{Levi form} of $\theta$ defined by $L_{\theta}(V,\overline{W}):=
-2id\theta(V\wedge \overline{W})$ on $\HH^M$ is definite we say that $\theta$ 
is
a \textit{pseudohermitian structure} (or \textit{contact form}), and that
$M$ is \textit{strictly pseudoconvex} with respect to $\theta$. A CR manifold 
with a fixed pseudohermitian structure is called a 
\textit{pseudohermitian manifold}. Given a pseudohermitian manifold 
$(M, \theta)$, by the nondegeneracy of the Levi form on $\HH^M$, 
$d\theta$ has precisely one null direction transverse to $\HH^M$ 
whence there is a uniquely defined real vector field $T$, called the 
\textit{characteristic vector field} of $\theta$, that satisfies
$$
T\hook d\theta = 0;\quad\theta(T) = 1,
$$
as in~\cite{lee86}.
A complex-valued $q$-form $\mu$ is said to be of type $(q,0)$ if 
$\overline{\HH^M}\hook\mu = 0$, and of type $(0,q)$ if $\HH^M\hook\mu = 0$.
An \textit{admissible coframe} on an open set of $M$ is a set of complex
 $(1,0)$-forms $\{\theta^{\alpha}\} = \{\theta^1,\dots,\theta^n\}$ that 
satisfy $\theta^\al(T) = 0$ and whose 
restriction to $\HH^M$ forms a basis for $(\HH^M)^\ast$.
 
The pseudohermitian manifolds
 we deal with in this note will be smooth compact CR manifolds of dimension
$2n+1$ embedded into a complex manifold $X$ of complex dimension $n+1$, 
with the natural CR structure 
$\HH^M = T_{1,0}X\cap\C TM$ and the pseudohermitian structure $\theta = 
(i/2)(\overline{\partial}\phi - 
\partial\phi)|_{M}$, where $\phi$ is a defining function for $M$. Here and 
throughout this note the pullback of forms via the inclusion of $M$ into $X$ 
is denoted by restriction $|_{M}$. The region $\{\phi < 0\}$
enclosed by $M$ will also be denoted by $X$. We shall assume that the Levi form
of $\theta$ is \textit{positive} definite.

\subsection{The ambient connection of Graham--Lee} There is a canonical 
connection on $M$ defined by 
Tanaka~\cite{tanaka} and Webster~\cite{webster} which, along with its curvature
 and torsion tensors, will be called the 
\textit{pseudohermitian connection}, respectively, 
\textit{pseudohermitian curvature}, \textit{pseudohermitian torsion}.
Now follows an exposition of the work of 
Graham--Lee~\cite{gl} in extending the pseudohermitian connection to a 
family of hypersurfaces that foliate a one-sided neighbourhood of $M$. 

Let $\Ub = U\cup M$ be a one-sided neighbourhood of $M$ where, for an 
arbitrary defining 
function $\phi$ for $M = M^0$, each 
level set $M^{\epsilon} := \{z: \phi(z) = \epsilon\}$ is strictly 
pseudoconvex
 with the natural CR structure 
$\HH^\epsilon := T_{1,0}\Ub\cap\C TM^\epsilon$ and pseudohermitian 
structure 
$(i/2)(\overline{\partial}\phi - \partial\phi)|_{M^\epsilon}$. Let 
$\HH\subset\C T\Ub$ denote the bundle whose fibre over each 
$M^{\epsilon}$ is $\HH^{\epsilon}$.  Setting $\vartheta := 
(i/2)(\overline{\partial}\phi - \partial\phi)$, the restriction of 
$i\partial\overline{\partial}\phi = d\vartheta$ to $\HH$ is positive definite. 
This means $\partial\overline{\partial}\phi$ has 
precisely one null direction transverse to $\HH$, whence there is a 
uniquely defined $(1,0)$ vector field $\xi$ that satisfies
\begin{equation}
\xi\perp_{\partial\overline{\partial}\phi}\HH;\quad\partial\phi(\xi) 
= 1.
\label{eq:xi}
\end{equation}

Define a  real-valued function $r := 2\ddbar\phi(\xi\wedge\overline{\xi})$ 
and call it
the \textit{transverse curvature} of $\phi$. Let 
$\{W_{\alpha}\} = \{W_1,\dots,W_n\}$ be any local frame for $\HH$. Since 
$\xi$ is transverse to $\HH$, the set of vector fields $\{W_{\alpha}, \xi\}$ 
is a local frame for $T_{1,0}\Ub$. The dual $(1,0)$ coframe is then 
of the form $\{\vartheta^\alpha, \partial\phi\}$ for some $(1,0)$-forms 
$\{\vartheta^{\alpha}\}$ that annihilate $\xi$. From the definitions of $\xi$ 
and $r$ we may 
write 
\begin{equation}
\ddbar\phi = h_{\al\bbar}\vartheta^{\alpha}\wedge\vartheta^{\bbar} 
+ r\partial\phi
\wedge\dbar\phi,
\label{eq:ddbarphi}
\end{equation}
for a positive definite Hermitian matrix of functions $h_{\al\bbar}$. We will
 use $h_{\al\bbar}$ and its inverse $h^{\al\bbar}$ to raise and lower 
indices in the usual way. Since 
\begin{equation}
d\vartheta = i\ddbar\phi = ih_{\al\bbar}\vartheta^{\alpha}\wedge
\vartheta^{\bbar} + 
rd\phi\wedge\vartheta,
\label{eq:dtheta}
\end{equation}
the Levi form is given by
$$
L_{\vartheta}(X,\overline{Y}) = -2id\vartheta(X\wedge\overline{Y}) = 
h_{\al\bbar}X^{\alpha}Y^{\overline{\beta}},
$$
for $X = X^{\alpha}W_{\alpha}, Y = Y^{\beta}W_{\beta}\in\HH$. Writing 
\begin{equation}
\label{eq:xi2}
\xi = N - \frac{i}{2}T
\end{equation} 
with $N$ and $T$ real, it follows that $T$ is tangential to each $M^\epsilon$
and that
$T\hook d\vartheta = 0$ and 
$\vartheta(T) = 1$. 
Thus on each $M^\epsilon$, $T$ is the characteristic vector field of 
$\vartheta$ and the restriction of 
$\{\vartheta^{\alpha}\}$ is an admissible coframe.  The 
pseudohermitian connection on each $M^\epsilon$, denoted here by 
$\nabla^{\epsilon}$, is
defined in~\cite{webster} by the relations
$$
\nabla^{\epsilon}W_{\al} = \omega_{\al}^{\phantom{\al}\beta}\otimes W_{\beta}
\:,
\quad\nabla^{\epsilon}T = 0\:;
$$
$$
d\vartheta^{\beta} = \vartheta^{\alpha}\wedge\omega_{\al}
^{\phantom{\al}\beta} + 
\vartheta\wedge\tau^{\beta}\:,
\quad \tau^{\beta} = A^{\beta}_{\phantom{\beta}\abar}\vartheta
^{\overline{\al}}\:,
\quad 
\omega_{\al\bbar} + \omega_{\bbar\al} = dh_{\al\bbar}\:,
$$
for a matrix of one-forms $\omega_{\al}^{\phantom{\al}\beta}$ and a 
symmetric matrix of 
functions $A_{\al\beta}$ on $M^{\epsilon}$. We 
extend 
$\omega_{\al}^{\phantom{\al}\beta}$ and $\tau^{\beta}$ uniquely to smooth 
one-forms on $\Ub$ by
 declaring $\omega_{\al}^{\phantom{\al}\beta}(N) = \tau^{\beta}(N) = 0$.

Graham--Lee proceed to define a connection, called the 
\textit{ambient connection}, on $\Ub$ that restricts to the 
pseudohermitian connection on each $M^\epsilon$.
\begin{prop}[{\cite[Proposition 1.1]{gl}}]
With $\Ub$ and $\phi$ as above, there exists a unique linear connection 
$\nabla$ on $\Ub$ with the properties:\\
\textnormal{(a)} For any vector fields $X$ and $Y$ on $\Ub$ that are tangent 
to some 
$M^\epsilon$, $\nabla_XY = \nabla^\epsilon_XY$, where $\nabla^\epsilon$ is 
the pseudohermitian connection on $M^\epsilon$.\\ 
\textnormal{(b)} $\nabla$ preserves $\HH, N, T $ and $L_{\vartheta}$.\\
\textnormal{(c)} Let $\{W_\alpha\}$ be any frame for $\HH$ and 
$\{\vartheta^\alpha, 
\partial\phi\}$ the $(1,0)$ coframe dual to $\{W_{\alpha}, \xi\}$. The 
connection 
one-forms  $\phi_{\al}^{\phantom{\al}\beta}$, defined by 
$\nabla W_{\al} = \phi_{\al}^{\phantom{\al}\beta}\otimes W_{\beta}$, satisfy 
the following structure equation:
\begin{equation}
d\vartheta^{\al} = \vartheta^{\beta}\wedge\phi_{\beta}^{\phantom{\beta}\al} - 
i\partial\phi\wedge\tau^{\al} + iW^\al rd\phi\wedge\vartheta 
+ \frac{1}{2}rd\phi
\wedge\vartheta^\alpha.
\label{eq:glstr}
\end{equation}
\label{prop:glstr}
\end{prop}
The components of successive covariant derivatives of a tensor with 
respect to $\nabla$ will appear as subscripts preceded by a comma, e.g., 
$A_{\alpha\beta,\:\overline{\gamma}\delta}\:$. For a scalar function the comma 
will be omitted. Of course covariant derivatives of a tensor in the 
direction $W_\alpha$ or $W_{\abar}$ can be interpreted as either referring to
 the pseudohermitian connection or the ambient connection. Covariant 
differentiation in the direction $W_0 := \xi$ shall be written in components
 as, e.g., $A_{\alpha\beta,0} = A_{\alpha\beta,N} - (i/2)A_{\alpha\beta,T}\:$.

\begin{prop}[{\cite[Proposition 1.2]{gl}}]
The curvature forms $\Omega_{\al}^{\phantom{\al}\beta} = 
d\phi_{\al}^{\phantom{\al}\beta} - \phi_{\al}^{\phantom{\al}\gamma}\wedge
\phi_{\gamma}^{\phantom{\gamma}\beta}$ are given by
\begin{eqnarray}
\label{eq:glcur}
\Omega_{\al}^{\phantom{\al}\beta} & = &
R^{\phantom{\al}\beta}
_{\al\phantom{\beta}\rho\overline{\gamma}}\vartheta^\rho\wedge\vartheta^
{\overline{\gamma}}
+
i\vartheta_\alpha\wedge\tau^\beta
-
i\tau_\alpha\wedge\vartheta^\beta\\
&&{}+
iA_{\alpha\gamma ,}^{\phantom{\alpha\gamma ,}\beta}\vartheta^\gamma\wedge
\dbar{\phi}
+
iA^{\beta}_{\phantom{\beta}\overline{\gamma},\alpha}\vartheta
^{\overline{\gamma}}
\wedge\partial{\phi}\nonumber\\
&&{}+
d\phi\wedge(r_\alpha\vartheta^\beta - r^\beta\vartheta_\alpha 
+ \frac{1}{2}\delta^
\beta_
\alpha r_\gamma\vartheta^\gamma - \frac{1}{2}\delta^\beta_
\alpha r_{\overline{\gamma}}\vartheta^{\overline{\gamma}})
\nonumber\\
&&{}+
d\phi\wedge\frac{i}{2}(r_{\al}^{\phantom{\al}\beta} + 
r^{\beta}_{\phantom{\beta}\alpha} + 2A_{\al\gamma}A^{\gamma\beta})\vartheta\:,
\nonumber
\end{eqnarray}
where $R^{\phantom{\al}\beta}_{\al\phantom{\beta}\rho\overline{\gamma}}$ 
agree with the components of the pseudohermitian curvature tensor.
\label{prop:glcur}
\end{prop}

For later use 
we define the real \textit{sublaplacian} on a function $f$ by
\begin{equation*}
\Delta_b f := -(f_\alpha^{\phantom{\alpha}\alpha} + 
f_{\bbar}^{\phantom{\bbar}\bbar}).
\end{equation*} We will also use a \textit{divergence formula} for 
pseudohermitian manifolds, which we state here for reference. For a 
$(1,0)$-form 
$\mu = \mu_\alpha\theta^\alpha$ on $M$, an application of Stokes' theorem to 
the $2n$-form $\theta\wedge\mu\wedge (d\theta)^{n-1}$ yields the formula
\begin{equation}
\label{eq:div}
\int_M\mu_{\alpha ,}^{\phantom{\alpha ,}\alpha}\ \theta\wedge (d\theta)^n = 0.
\end{equation}

\section{The Einstein equation}
\label{sec:einstein}

\subsection{Choice of defining function and frame} 
\label{subsec:frame}
Given any defining function 
$\phi$ for $M$ we may identify (a possible shrinkage of) $\Ub$ 
with $M \times (\delta,0]$, $\delta<0$, with $\phi$ as a coordinate on the 
second factor. Indeed, identify $\Ub$ with $M \times (\delta,0]$ 
by the map $F$ that 
takes $(p,t)$ to the point in $\Ub$ obtained by following the unit-speed 
integral curve $\gamma^p(t)$ of the vector field $N$, 
emanating from a point $p$ on $M$, for time $t$. The remaining assertion is 
equivalent to the statement that $\gamma^p(t)\in\{\phi = t\} = M^t$, which  
follows from the fact that $N\phi= 1$.

Specifying a defining 
function $\phi$ for $M$ determines a pseudohermitian structure 
$\theta$ for $M$ such that 
$\theta = (i/2)(\dbar\phi - \partial\phi)|_M$. In the other direction
though, the boundary equation 
\begin{equation} 
\frac{i}{2}(\dbar\phi - \partial\phi)|_M = \theta,
\label{eq:neumann}
\end{equation}
for a specified pseudohermitian structure $\theta$, clearly determines
$\phi$ only mod $O(\phi^2)$. (Throughout this note $O(\phi^k)$ will denote 
functions smoothly divisible by $\phi^k$.) However if there is a certain form 
of complete 
K\"ahler metric on $U$, 
then $\phi$ can be uniquely determined as follows. For any fixed sufficiently 
smooth defining function $\rho$ with 
$\rho < 0$ in $U$, consider the K\"ahler metric $g_+$ with K\"ahler form
\begin{equation}
\omega = \partial\overline{\partial}\log\left(-\frac{1}{\phicy}\right)
       = -\frac{\ddbar\phicy}{\phicy} 
+ \frac{\partial\phicy\wedge\dbar\phicy}{\phicy^2}.
\label{eq:kahlerform}
\end{equation}
If we write our desired defining function $\phi$ as $\phi = e^{2f}\rho$, for a 
function $f$, then using (\ref{eq:ddbarphi}),
\begin{equation}
\label{eq:metric}
\omega = -\frac{1}{\phi}h_{\al\bbar}\vartheta^{\alpha}\wedge
\vartheta^{\bbar} + 
\frac{1-r\phi}{\phi^2}\partial\phi\wedge\dbar\phi + 2\ddbar f.
\end{equation}
This formula implies that
$$
|\partial(\log(-\phi))|_{g_+}^2 = 1 + O(\phi).
$$ 
(Our convention is that a K\"ahler metric 
$g_{i\overline{j}}\vartheta^i\vartheta^{\overline{j}}$ has 
K\"ahler form $g_{i\overline{j}}\vartheta^i\wedge\vartheta^{\overline{j}}$.)
Then in analogy with the conformal case, we choose a special
 defining function $\phi$ according to the following lemma, whose proof is 
similar to that of Lemma 2.1 in~\cite{graham}.
\begin{lem}
A choice of pseudohermitian structure $\theta$ on $M$ determines a unique 
defining function $\phi$ in a neighbourhood of $M$ such that 
$$
\frac{i}{2}(\dbar\phi - \partial\phi)|_M = \theta\ ;\quad
|\partial(\log(-\phi))|_{g_+}^2 = 1.
$$
\label{lem:special}  
\end{lem}
\begin{proof}
As before write $\phi = e^{2f}\phicy$ for
a function $f$ to be determined. The boundary value of $f$ is determined by 
the condition $(i/2)(\dbar\phi - \partial\phi)|_M = \theta$. Next
$$
|\partial(\log(-\phi))|^2_{g_+} = |\partial(\log(-\rho))|^2_{g_+} 
+ \frac{4}{\rho}
\Re \langle\partial\rho,\partial f\rangle_{g_+} + 4|\partial f|^2_{g_+}.
$$
If $N_\rho$ denotes the real part of the vector field that solves (\ref{eq:xi})
 with the defining function $\rho$, then from (\ref{eq:kahlerform}) and 
(\ref{eq:ddbarphi}) (for the defining function $\rho$) the equation 
above becomes 
$$
|\partial(\log(-\phi))|^2_{g_+} = \frac{1}{\rho^2}|\partial\rho|^2_{g_+} 
+ \frac{4\rho}{1-r\rho}N_\rho f +  4|\partial f|^2_{g_+}.
$$
The condition $|\partial(\log(-\phi))|^2_{g_+} = 1$ is equivalent to
$$
N_\rho f +\frac{(1 - r\rho)}{\phicy}|\partial f|^2_{g_+} 
= \frac{(1 - \frac{1}{\phicy^2}|\partial\phicy|^2_{g_+})(1-r\rho)}{4\phicy}\:.
$$
The numerator of the right-hand side is $O(\rho)$. Moreover 
the second term on the left-hand side is $O(1)$ with the coefficient of 
$(N_\rho f)^2$ being $O(\rho)$. Therefore we have a 
noncharacteristic first-order PDE for
 $f$ and a unique solution near $M$ with the prescribed boundary condition.
\end{proof}
Hereafter $\phi$ will denote the special defining function produced by 
Lemma~\ref{lem:special}. 

We want to find 
$(1,0)$-forms $\{\tilde{\vartheta}^\alpha\}$ so that, 
keeping Lemma~\ref{lem:special} in mind, the metric
$g_+$ has the diagonal K\"ahler form
\begin{equation}
\omega = -\frac{1}{\phi}\hhat_{\al\bbar}\varthetat^{\alpha}\wedge
\varthetat^{\bbar} + 
\frac{1}{\phi^2}\partial\phi\wedge\dbar\phi,
\label{eq:diagmetric} 
\end{equation}
where $\hhat_{\al\bbar}$ is a positive definite Hermitian matrix of functions,
 parameterised by $\phi$. 
 To construct the desired forms, first take the unique vector field $\xit$ 
of type $(1,0)$ that satisfies 
$$
\xit\perp_{g_+}\HH;\quad\partial\phi(\xit) = 1.
$$
Thus $\{W_{\alpha}, \xit\}$ is a local frame for $T_{1,0}U$. The dual 
$(1,0)$ coframe is then of the form $\{\tilde{\vartheta}^\alpha, 
\partial\phi\}$ 
for some $(1,0)$ forms $\{\tilde{\vartheta}^{\alpha}\}$ on 
$U$ that annihilate $\xit$. We may write 
$$
\tilde{\vartheta}^\alpha = \vartheta^\alpha + a^\alpha\partial\phi
$$
and
$$
\xit = \xi - a^\alpha W_\alpha 
$$
for some uniquely determined functions $a^\alpha$ on $U$.

\begin{lem}
The functions $a^\alpha$ extend continuously to $M$ with $a^\alpha|_M = 0$.
Consequently, the vector field $\xit$, the forms $\varthetat^\alpha$ and the
tensor $\hhat_{\alpha\bbar}$ have continuous extensions to $M$ with
$\xit|_M = \xi|_M, \varthetat^\alpha|_M = \vartheta^\alpha|_M$ and 
$\hhat_{\al\bbar}|_M = h_{\al\bbar}|_M$. 
\end{lem}
\begin{proof}
Using (\ref{eq:diagmetric}) and (\ref{eq:metric}),
\begin{equation*}
\begin{split}
0 = 2\omega(\xit\wedge W_{\bbar}) &= 2\omega(\xi\wedge W_{\bbar}) 
- 2a^\alpha\omega(W_\alpha\wedge W_{\bbar})\\
&= 4\ddbar f(\xi\wedge W_{\bbar}) - 4a^\alpha \ddbar f (W_\alpha\wedge 
W_{\bbar}) + \phi^{-1} a^\alpha h_{\alpha\bbar}.
\end{split}
\end{equation*}
This implies that $a^\alpha$ is 
$O(\phi)$, whence has a continuous extension to the boundary by setting 
$a^\alpha|_M = 0$. The assertions for $\xit, \varthetat^\alpha$ and 
$\hhat_{\alpha\bbar}$ follow immediately.
\end{proof}
Defining the vector field $\eta = \eta^\alpha W_\alpha$ for $\eta^\alpha := 
\phi^{-1}a^\alpha$, we may rewrite the equations for $\varthetat^\alpha$ and 
$\xit$  as
\begin{equation}
\tilde{\vartheta}^\alpha = \vartheta^\alpha + \phi \eta^\alpha\partial\phi
\label{eq:thetat}
\end{equation}
and
\begin{equation}
\label{eq:eta}
\xit = \xi - \phi \eta
\end{equation}
in $\Ub$.

\begin{lem} 
For the special defining function chosen according to Lemma~\ref{lem:special},
the transverse curvature function $r$ vanishes on $M$.
\end{lem}
\begin{proof}
Using (\ref{eq:diagmetric}) and (\ref{eq:metric}),
\begin{equation*}
\frac{1}{\phi^2} = 2\omega(\xit\wedge\overline{\xit}) = 
-\frac{2}{\phi}\ddbar\phi(\xit\wedge\overline{\xit}) + \frac{1}{\phi^2} 
+ 4\ddbar f
(\xit\wedge\overline{\xit}).
\end{equation*}
Hence $\ddbar\phi(\xit\wedge\overline{\xit})$ is $O(\phi)$ from which it 
follows, by (\ref{eq:eta}) and the definition of $r$, that $r$ is $O(\phi)$.
\end{proof}

\subsection{Levi-Civita connection and Ricci form} We use the notation $D$ 
for the 
Levi-Civita connection of $(X,g_+)$. We shall compute the connection matrix 
$(\psi_j^{\phantom{j}k})$ of $D$ with respect to the coframe 
$\{\tilde{\vartheta}^j\}$, where $\tilde{\vartheta}^ 0 := \partial\phi$. 
Denote by $g_{j\kb}$ the components of $g_+$ with respect to this 
coframe, with 
$$
(g_{j\kb}) = 
\left(\begin{array}{cc}
\phi^{-2} & 0 \\
0 & 
-\phi^{-1}\hhat_{\alpha\bbar}
\end{array}\right),
$$
$$
(\ginv^{j\kb}) = 
\left(\begin{array}{cc}
\phi^{2} & 0 \\
0 & 
-\phi\hhatinv^{\alpha\bbar}
\end{array}\right).
$$
The matrix
$(\psi_j^{\phantom{j}k})$ is uniquely determined by the structure equations
\begin{equation}
d\varthetat^j = \varthetat^k\wedge\psi_k^{\phantom{k}j};
\label{eq:struct}
\end{equation}
\begin{equation}
\psi_{j\kb} + \psi_{\kb j} = dg_{j\kb}.
\label{eq:compat}
\end{equation}
 Substituting  
$\psi_k^{\phantom{k}j} = \psi_{k\phantom{j}l}^{\phantom{k}j}\varthetat^l 
+ \psi_{k\phantom{j}\lb}^{\phantom{k}j}\varthetat^{\lb}$ 
into 
(\ref{eq:struct}) and using (\ref{eq:ddbarphi}), (\ref{eq:glstr}) and 
(\ref{eq:thetat}) will uniquely determine all of the 
$\psi_{k\phantom{j}\lb}^{\phantom{k}j}$. The $\psi_{k\phantom{j}l}
^{\phantom{k}j}$ can then
be uniquely determined from (\ref{eq:compat}). The result of performing these
computations is the following:
\begin{prop}
The Levi-Civita connection matrix of $g_+$ with respect to the coframe
$\{\tilde{\vartheta}^j\}$ is
\label{prop:connection}
$$
(\psi_j^{\phantom{j}k}) = 
\left(\begin{array}{cc}
\psi_ 0^{\phantom{ 0} 0} & \psi_\beta^{\phantom{\beta} 0} \\
\psi_ 0^{\phantom{ 0}\alpha} & \psi_\beta^{\phantom{\beta}\alpha}
\end{array}\right),
$$
where 
\begin{eqnarray*}
\psi_ 0^{\phantom{ 0} 0} & = &\left(-2\phi^{-1} + r 
+ \phi^2 |\eta|^2  \right)\partial\phi\\
&& + \left(-\phi \eta_{\alpha}\right)\varthetat^\alpha\\
&& + \left(- r - \phi^2 |\eta|^2  \right)\dbar\phi\\
&& + \left(\phi \eta_{\bbar}  \right)\varthetat^{\bbar}\:;\\
\psi_\beta^{\phantom{\beta} 0} & = 
&\left(-\phi \eta_{\beta}  \right)\partial\phi\\
&& + \left(\phi\hhat_{\beta\gbar}
\left(iA^{\gbar}_{\phantom{\gbar}\alpha} - 
\phi W_\alpha \eta^{\gbar}
- \phi \eta^{\delbar}\overline{\phi_{\delta\phantom{\gamma}\abar}
^{\phantom{\delta}\gamma}}
+ \phi^2 \eta^{\gbar}\eta_{\alpha}
  \right)   \right)\varthetat^\alpha\\
&& + \left(\phi \eta_{\beta}  \right)\dbar\phi\\
&& + \left(-h_{\beta\gbar}  \right)\varthetat^{\gbar}\:;\\
\psi_ 0^{\phantom{ 0}\alpha} & = 
&\left(\hhatinv^{\alpha\gbar}\eta_{\gbar}  
\right)\partial\phi\\
&& + \left(-\phi^{-1}\hhatinv^{\alpha\gbar}h_{\delta\gbar}\right)\varthetat^
\delta\\
&& + \Big(-r^\alpha - \eta^\alpha - \phi\overline{\xi}\eta^\alpha
-\phi \eta^\alpha r + \phi \eta^{\bbar}(iA^\alpha_{\phantom{\alpha}\bbar} 
+ \phi W_{\bbar} \eta^\alpha)\\
&&\quad\ \  - \phi \eta^\beta
(\phi_{\beta\phantom{\alpha}\zerob}^{\phantom{\beta}\alpha} 
- \delta^\alpha_\beta\tfrac{1}{2}r 
- \phi \eta^{\gbar}\phi_{\beta\phantom{\alpha}\gbar}^{\phantom{\beta}\alpha}
+ \phi^2 \eta^\alpha \eta_{\beta}  )  \Big)\dbar\phi\\
&& + \left(-iA^\alpha_{\phantom{\alpha}\bbar} - 
\phi W_{\bbar}\eta^\alpha
- \phi \eta^{\gamma}\phi_{\gamma\phantom{\alpha}\bbar}
^{\phantom{\gamma}\alpha} 
+ \phi^2 \eta^\alpha \eta_{\bbar} \right)\varthetat^{\bbar}\:;\\
\psi_\beta^{\phantom{\beta}\alpha} & = 
&\Big(\phi\hhatinv^{\alpha\gbar}\left(
\xi (\phi^{-1}\hhat_{\beta\gbar}) - \phi^{-1}\hhat_{\beta\delbar}
(\overline{\phi_{\gamma\phantom{\delta}\zerob}^{\phantom{\gamma}\delta}} - 
\delta_{\gbar}^{\delbar}\tfrac{1}{2}r)\right) \\
&&\ \  - \phi \eta^\gamma \hhatinv^{\alpha\delbar}
\left(W_\gamma\hhat_{\beta\delbar} - \hhat_{\beta\mbar}
(\overline{\phi_{\delta\phantom{\mu}\gbar}^{\phantom{\delta}\mu}} 
- \phi \eta^{\mbar}
h_{\gamma\delbar})\right)\Big)\partial\phi\\
&& + \left(\hhatinv^{\alpha\gbar}\left(W_\delta\hhat_{\beta\gbar}  -
\hhat_{\beta\mbar}
(\overline{\phi_{\gamma\phantom{\mu}\delbar}^{\phantom{\gamma}\mu}}
- \phi \eta^{\mbar}h_{\delta\gbar}) \right)  
\right)\varthetat^\delta\\
&& + \left(\phi_{\beta\phantom{\alpha}\zerob}^{\phantom{\beta}\alpha}
- \delta_\beta^\alpha\tfrac{1}{2}r 
- \phi \eta^{\gbar}\phi_{\beta\phantom{\alpha}\gbar}^{\phantom{\beta}\alpha} 
+ \phi^2 \eta^\alpha \eta_{\beta}  \right)\dbar\phi\\
&& + \left(\phi_{\beta\phantom{\alpha}\gbar}^{\phantom{\beta}\alpha} 
- \phi \eta^\alpha h_{\beta\gbar}  \right)\varthetat^{\gbar}\:.\\
\end{eqnarray*}
Here we have written the ambient connection forms as
$$
\phi_\beta^{\phantom{\beta}\alpha} =
\phi_{\beta\phantom{\alpha} 0}^{\phantom{\beta}\alpha}\partial\phi
+\phi_{\beta\phantom{\alpha}\delta}^{\phantom{\beta}\alpha}\vartheta^\delta
+ \phi_{\beta\phantom{\alpha}\zerob}^{\phantom{\beta}\alpha}\dbar\phi
+ \phi_{\beta\phantom{\alpha}\gbar}^{\phantom{\beta}\alpha}\vartheta^{\gbar}.
$$
\label{prop:levicivita}
\end{prop}

\begin{prop}
The Ricci form $\rho$ of $g_+$ is
$\rho = \rho_{j\kb}\:\varthetat^{j}\wedge\varthetat^{\kb}$ 
where 
\begin{eqnarray*}
\rho_{ 0\zerob} &= &-(n + 2)\left(\frac{1}{\phi^2} + Nr -
\frac{1}{\phi}r\right) \\
&&  - (n + 2)r^2 - r\hhatinv^{\delta\gbar}\hhat_{\delta\gbar ,\: 0}
-  r^\alpha\hhatinv^{\delta\gbar}\hhat_{\delta\gbar ,\:\alpha}
-  (\hhatinv^{\delta\gbar}\hhat_{\delta\gbar ,\: 0})_{,\:\zerob}\\
&&  + \frac{1}{2}(\Delta_b r - 2|A|^2)\\
&&  - \phi^2 \eta^\alpha \eta^{\bbar}\rho_{\alpha\bbar};\\
\rho_{\alpha\zerob} &= &-\phi\rho_{\alpha\bbar}\eta^{\bbar}-
(\hhatinv^{\delta\gbar}\hhat_{\delta\gbar ,\:\alpha})_{,\: \zerob} +
iA_{\beta\alpha,}^{\phantom{\beta\alpha,}\beta}\\ 
&&  -(n+2)r_\alpha 
-\frac{1}{2}r\hhatinv^{\delta\gbar}\hhat_{\delta\gbar ,\:\alpha} \:;\\
\rho_{\alpha\bbar} &= &(n + 2)\frac{1}{\phi}h_{\al\bbar} \\
&& - (n+2)rh_{\al\bbar} - 
(\hhatinv^{\delta\gbar}\hhat_{\delta\gbar ,\: 0})h_{\al\bbar} -
(\hhatinv^{\delta\gbar}\hhat_{\delta\gbar ,\:\alpha})_{,\:\bbar} + 
\Ric_{\alpha\bbar}\:,\\
\end{eqnarray*}
and $\Ric_{\alpha\bbar}$ denotes the ambient Ricci tensor. 
\label{prop:ricci}
\end{prop}
\begin{proof}
Recall that the Ricci form of a K\"ahler metric is defined to be 
$d\psi_j^{\phantom{j}j}$. Using Proposition~\ref{prop:levicivita}, rewriting 
in terms of covariant derivatives and using the commutation relations 
in Lemma~\ref{lem:comm2} below
gives the result. (One can also use Lemma~\ref{lem:comm2} to check 
that Ricci form as given by Proposition
~\ref{prop:ricci} is in fact Hermitian symmetric.)
\end{proof}

\subsection{Einstein equation} As mentioned in the Introduction, 
Fefferman~\cite{fef76} solved the approximate Einstein equation
$$
\rho = -(n+2)\omega + O(\phi^{n+1}).
$$
We shall use this equation to determine our desired power series expansion for
$g_+$; in the following subsection we relate this expansion to Fefferman's 
solution. Now the $0\zerob$ component of this equation is
\begin{equation}
\begin{split}
&(n + 2)\phi N r - (n + 2)r\\
&+\phi\Big((n + 2)r^2 + r\hhatinv^{\delta\gbar}\hhat_{\delta\gbar ,\: 0}
+  r^\alpha\hhatinv^{\delta\gbar}\hhat_{\delta\gbar ,\:\alpha}
+  (\hhatinv^{\delta\gbar}\hhat_{\delta\gbar ,\: 0})_{,\:\zerob}\\
&\quad\quad\   - \frac{1}{2}(\Delta_b r - 2|A|^2)\Big)\\
&  + \phi^2(n + 2)\eta^\alpha \eta^{\bbar}\hhat_{\alpha\bbar}=0\:;
\label{eq:einsteinii}
\end{split}
\end{equation}
the $\alpha\zerob$ component is
\begin{equation}
\begin{split}
(n+2)\hhat_{\alpha\bbar}\eta^{\bbar} +
(\hhatinv^{\delta\gbar}\hhat_{\delta\gbar ,\:\alpha})_{,\: \zerob} -
iA_{\beta\alpha,\:}^{\phantom{\beta\alpha\:}\beta}
 +(n+2)r_\alpha +
\frac{r}{2}\hhatinv^{\delta\gbar}\hhat_{\delta\gbar ,\:\alpha} = 0\:;
\label{eq:einsteinai}
\end{split}
\end{equation}
and the $\alpha\bbar$ component is 
\begin{equation}
\begin{split}
&(n + 2)(\hhat_{\al\bbar} - h_{\al\bbar})\\
& + \phi\Big((n+2)rh_{\al\bbar} + 
(\hhatinv^{\delta\gbar}\hhat_{\delta\gbar ,\: 0})h_{\al\bbar} +
(\hhatinv^{\delta\gbar}\hhat_{\delta\gbar ,\:\alpha})_{,\:\bbar} - 
\Ric_{\alpha\bbar}\Big)  = 0\:.
\label{eq:einsteinab}
\end{split}
\end{equation}
Substituting the formula for $\eta^\alpha$ from (\ref{eq:einsteinai}) into 
(\ref{eq:einsteinii}) gives
\begin{equation}
\begin{split}
&(n + 2)\phi N r - (n + 2)r\\
&+\phi\Big((n + 2)r^2 + r\hhatinv^{\delta\gbar}\hhat_{\delta\gbar ,\: 0}
+  r^\alpha\hhatinv^{\delta\gbar}\hhat_{\delta\gbar ,\:\alpha}
+  (\hhatinv^{\delta\gbar}\hhat_{\delta\gbar ,\: 0})_{,\:\zerob}\\
&  \quad\quad\ - \frac{1}{2}(\Delta_b r - 2|A|^2)\Big)\\
&  + \phi^2\frac{\hhatinv^{\alpha\bbar}}{n+2}\Big(
(\hhatinv^{\delta\gbar}\hhat_{\delta\gbar ,\:\bbar})_{,\: 0} +
iA_{\gbar\bbar,\:}^{\phantom{\gbar\bbar\:}\gbar}
 +(n+2)r_{\bbar} + 
\frac{1}{2}r\hhatinv^{\delta\gbar}\hhat_{\delta\gbar ,\:\bbar}\Big)\\
&\times\Big((\hhatinv^{\delta\gbar}\hhat_{\delta\gbar ,\:\alpha})_{,\: \zerob}
 -iA_{\beta\alpha,\:}^{\phantom{\beta\alpha\:}\beta}
 +(n+2)r_\alpha + 
\frac{1}{2}r\hhatinv^{\delta\gbar}\hhat_{\delta\gbar ,\:\alpha}\Big)=0\:,
\label{eq:einsteinii2}
\end{split}
\end{equation} 
and then our task is to use the system (\ref{eq:einsteinab},
~\ref{eq:einsteinii2}) of $n(n+1)/2 + 1$ equations to obtain the Taylor series 
coefficients of the functions $\hhat_{\alpha\bbar}$. As all the terms here
transform tensorially, we can 
take ambient connection covariant derivatives as often as we please. 
To determine the components of the tensors 
$\nabla_N\hhat|_M, \nabla_N^2\hhat|_M,$ etc., we use the following:
\begin{lem}
\label{lem:spframe}
Given a point $p\in \Ub$ there exists a frame 
$\{W_\alpha\}$ for $\HH$ in a neighbourhood of $p$ with respect to which 
$\phi_{\al\phantom{\al} N}^{\phantom{\al}\beta} := 
\phi_{\al}^{\phantom{\al}\beta}(N)\equiv 0$ near $p$.
\end{lem}
\begin{proof}
Let $\{\Wh_\alpha\}$ be any frame for $\HH$ defined near 
$p\in M^\epsilon$. 
Let $\{\varthetah^{\alpha}, \partial\phi\}$ be the dual coframe 
to $\{\Wh_\alpha, \xi\}$.
Set  $\vartheta^\alpha = a^\alpha_\beta\varthetah^\beta, 
W_\alpha = (a^{-1})^\beta_\alpha \Wh_\beta$ 
for a matrix of functions $(a^\alpha_\beta)$ uniquely determined as follows. 
Substituting
$a^\alpha_\beta\varthetah^\beta$ in place of $\vartheta^\alpha$ in the 
structure
equation (\ref{eq:glstr}) and examining only the $d\phi\wedge\varthetah^\beta$ 
component yields 
$$
Na^\alpha_\beta - 
a^\alpha_\gamma\phih_{\beta\phantom{\gamma} N}^{\phantom{\beta}\gamma} = 
- a^\delta_\beta\phi_{\delta\phantom{\alpha}N}^{\phantom{\delta}\al},
$$
where $\phih_{\al}^{\phantom{\al}\beta}$ and 
$\phi_{\al}^{\phantom{\al}\beta}$ are the ambient connection forms with 
respect to $\{\Wh_\alpha, \xi\}$ and 
$\{W_\alpha, \xi\}$ respectively. Set $(a^\alpha_\beta)$ to be the unique local
solution to the initial-value problem
\begin{equation}
N a^\alpha_\beta - 
a^\alpha_\gamma\phih_{\beta\phantom{\gamma} N}^{\phantom{\beta}\gamma} 
= 0\:;\quad
(a^\alpha_\beta)|_{M^\epsilon} = ((\overset{\circ}{a})^\alpha_\beta),
\label{eq:ivp2}
\end{equation}
for an arbitrarily prescribed matrix $((\overset{\circ}{a})^\alpha_\beta)$.
It is clear then that $\phi_{\delta\phantom{\alpha} N}
^{\phantom{\delta}\al}\equiv 0$ near
$p$. 
\end{proof}

Using a special frame produced by this lemma, 
ambient covariant differentiation in the direction 
$N$ is simply the ordinary directional derivative 
$N$. So the components of the tensors 
$\nabla_N\hhat|_M, \nabla_N^2\hhat|_M,$ etc., with respect to this
special frame are nothing more than
$N\hhat_{\alpha\bbar}, N^2\hhat_{\alpha\bbar},$ etc.. Moreover,
$h_{\alpha\bbar}$ becomes constant in the $N$ direction.

Before dealing with the system (\ref{eq:einsteinab},~\ref{eq:einsteinii2}) we 
gather some useful identities.
The proofs of the following three lemmas use the familiar result 
(see, e.g.,~\cite{lee88}) that there is a frame for $\HH$ in which the ambient 
connection forms vanish at a point. 

\begin{lem} The second covariant derivatives of a tensor of the form
$t = t_{\delta\gbar}$ satisfy the commutation relations
\begin{eqnarray}
t_{\delta\gbar ,\: N T} - t_{\delta\gbar ,\:T N} &=&
ir^\alpha t_{\delta\gbar,\:\alpha} - ir^{\bbar} t_{\delta\gbar,\:\bbar}
+ rt_{\delta\gbar,\: T}\\
&&+ \frac{i}{2}(r_{\delta}^{\phantom{\delta}\mu} +
r^{\mu}_{\phantom{\mu}\delta} + 2A_{\delta\alpha}A^{\alpha\mu})t_{\mu\gbar}
\nonumber\\
&&- \frac{i}{2}(r_{\gbar}^{\phantom{\gbar}\mbar} +
r^{\mbar}_{\phantom{\mbar}\gbar} + 2A_{\gbar\bbar}A^{\bbar\mbar})
t_{\delta\mbar}\nonumber;\\
t_{\delta\gbar ,\:N\alpha} - t_{\delta\gbar ,\:\alpha N} &=&
\frac{1}{2}r t_{\delta\gbar ,\:\alpha}
+ \frac{i}{2}A^{\bbar}_{\phantom{\bbar}\alpha} t_{\delta\gbar ,\:\bbar}\\
&&+ \frac{i}{2}A^{\mbar}_{\phantom{\mu}\alpha,\:\gbar} t_{\delta\mbar}
- \frac{i}{2}A_{\delta\alpha,\:}^{\phantom{\delta\alpha,\:}\mu} t_{\mu\gbar}
\nonumber\\
&&- t_{\delta\mbar}r^{\mbar} h_{\alpha\gbar} 
+ t_{\alpha\gbar}r_{\delta}\nonumber;\\
t_{\delta\gbar ,\: T\alpha} - t_{\delta\gbar ,\:\alpha T} &=&
A^{\bbar}_{\phantom{\bbar}\alpha} t_{\delta\gbar ,\:\bbar}
+ A^{\mbar}_{\phantom{\mbar}\alpha,\:\gbar} t_{\delta\mbar}
- A_{\delta\alpha,\:}^{\phantom{\delta\alpha,\:}\mu} t_{\mu\gbar};\\
t_{\delta\gbar ,\:\alpha\mu} - t_{\delta\gbar ,\:\mu\alpha} &=&
ih_{\delta\bbar}(A^{\bbar}_{\phantom{\bbar}\mu} t_{\alpha\gbar} -
A^{\bbar}_{\phantom{\bbar}\alpha} t_{\mu\gbar})\\
&&- it_{\delta\bbar}(A^{\bbar}_{\phantom{\bbar}\mu} h_{\alpha\gbar} -
A^{\bbar}_{\phantom{\bbar}\alpha} h_{\mu\gbar});\nonumber\\
t_{\delta\gbar ,\:\alpha\bbar} - t_{\delta\gbar ,\:\bbar\alpha} &=&
ih_{\alpha\bbar}t_{\delta\gbar ,\: T} + 
t_{\mu\gbar}R_{\delta\phantom{\mu}\alpha\bbar}^{\phantom{\delta}\mu}
- t_{\delta\mbar}R_{\gbar\phantom{\mbar}\bbar\alpha}^{\phantom{\gbar}\mbar}.
\label{eq:abbar}
\end{eqnarray}
\label{lem:comm2}
\end{lem}
\begin{proof}
Take the exterior derivative $d$ of the identity
$$
dt_{\delta\gbar} = t_{\delta\gbar,\:N}d\phi
+ t_{\delta\gbar,\: T}\vartheta
+ t_{\delta\gbar,\:\alpha}\vartheta^\alpha 
+ t_{\delta\gbar,\:\bbar}\vartheta^{\bbar}
+ t_{\mu\gbar}\phi_\delta^{\phantom{\delta}\mu}
+ t_{\delta\mbar}\phi_{\gbar}^{\phantom{\gbar}\mbar},
$$
use (\ref{eq:glstr}) and (\ref{eq:glcur}),
and compare like terms at a point where the ambient connection forms 
$\phi_\alpha^{\phantom{\alpha}\beta}$ vanish.
\end{proof}
\begin{lem} The third covariant derivatives of a tensor of the form
$t = t_{\delta\gbar}$ satisfy the commutation relations
\begin{eqnarray}
\quad\quad t_{\delta\gbar ,\:\alpha N\bbar} 
- t_{\delta\gbar ,\:\alpha\bbar N} &=&
\frac{1}{2}r t_{\delta\gbar ,\:\alpha\bbar}
- \frac{i}{2}A^\mu_{\phantom{\mu}\bbar} t_{\delta\gbar ,\:\alpha\mu}\\
&&- \frac{i}{2}A^\mu_{\phantom{\mu}\bbar,\:\delta} t_{\mu\gbar,\:\alpha}
+ \frac{i}{2}A_{\gbar\bbar,\:}^{\phantom{\gbar\bbar,\:}\mbar} 
t_{\delta\mbar,\:\alpha}
\nonumber\\
&&- t_{\mu\gbar,\:\alpha} r^\mu h_{\delta\bbar} 
+ t_{\delta\bbar,\:\alpha}r_{\gbar}\nonumber\\
&& -\frac{i}{2}A^\mu_{\phantom{\mu}\bbar,\:\alpha} t_{\delta\gbar,\:\mu}
- t_{\delta\gbar,\:\mu} r^\mu h_{\alpha\bbar}
- \frac{1}{2}t_{\delta\gbar,\:\alpha}r_{\bbar}\nonumber.
\end{eqnarray}
\label{lem:comm3}
\end{lem}
\begin{proof}
Similar to that of Lemma~\ref{lem:comm2}.
\end{proof}
\begin{rem}
If we take the tensor $t$ to be $\hhat$, because 
$\hhat_{\alpha\bbar}|_M = h_{\alpha\beta}|_M$, 
$\nabla h\equiv 0$ and $r|_M\equiv 0$, the right-hand sides of all the 
equations in 
Lemmas~\ref{lem:comm2} and~\ref{lem:comm3} vanish on $M$.
\end{rem}

\begin{lem} The ambient Ricci and torsion tensors satisfy
\begin{eqnarray}
\Ric_{\delta\gbar,\:N} &=& \frac{i}{2}(
A_{\alpha\delta,\phantom{\alpha}\gbar}^{\phantom{\alpha\delta,}\alpha}
- A_{\bbar\gbar,\phantom{\bbar}\delta}^{\phantom{\bbar\gbar,}\bbar})
- (\frac{n}{2} + 1)(r_{\delta\gbar} + r_{\gbar\delta})\label{eq:ricci}\\
&&+ \frac{1}{2}(\Delta_b r - 
2|A|^2)h_{\delta\gbar} - r\Ric_{\delta\gbar};\nonumber\\
A_{\alpha\beta,\:N} &=& -rA_{\alpha\beta} + ir_{\alpha\beta} + \frac{i}{2}
A_{\alpha\beta,\:T}.
\end{eqnarray}
\label{lem:bianchi}
\end{lem}
\begin{proof}
Work in a frame where the ambient connection forms vanish at a point. Then 
 the first identity follows by taking the exterior 
derivative $d$ of the curvature equation (\ref{eq:glcur}) and the second 
identity is obtained by taking $d$ of the structure equation (\ref{eq:glstr}).
\end{proof}

Now applying $\nabla_{N}$ to (\ref{eq:einsteinab}), evaluating on $M$ and
using a prime to denote covariant differentiation in the direction $N$ gives
$$
(n+2)\hhat_{\alpha\bbar}^\prime + h^{\delta\gbar}\hhat_{\delta\gbar}^\prime
h_{\alpha\bbar} - \Ric_{\alpha\bbar} = 0.
$$
Taking trace,
$$
h^{\alpha\bbar}\hhat_{\alpha\bbar}^\prime\Big|_M = \frac{1}{2(n+1)}\Scal,
$$
where $\Scal$ is the pseudohermitian scalar curvature. Substituting this into 
the previous equation yields
\begin{equation}
\label{eq:hhatprime}
\hhat_{\alpha\bbar}^\prime|_M = \frac{1}{n+2}\left(\Ric_{\alpha\bbar}
-\frac{1}{2(n+1)}\Scal h_{\alpha\bbar} \right).
\end{equation}
Remark that this is the so-called \textit{conformal Ricci tensor}; an 
analogous tensor appears in the conformal setting also (\cite{graham}).

In order to determine $\hhat^{\prime\prime}_{\alpha\bbar}|_M$ we first 
calculate $h^{\alpha\bbar}\hhat^{\prime\prime}_{\alpha\bbar}|_M$ in terms 
of $r^\prime|_M$ and pseudohermitian invariants of $(M,\theta)$ by applying
$\nabla_N$ to (\ref{eq:einsteinii2}) and evaluating on $M$. Applying 
$\nabla^2_N$ to (\ref{eq:einsteinab}), taking trace and inserting the 
just-computed value for $h^{\alpha\bbar}\hhat^{\prime\prime}_{\alpha\bbar}|_M$
determines $r^\prime|_M$ and hence $\hhat^{\prime\prime}_{\alpha\bbar}|_M$.

In general, applying $\nabla_N^k$ to (\ref{eq:einsteinab}), taking 
trace and evaluating on $M$ yields
\begin{equation}
h^{\delta\gbar}\hhat_{\delta\gbar}^{(k)}|_M =
\frac{-k n(n+2)r^{(k-1)}|_M}
{kn + n + 2} \: + \: \textnormal{(terms involving $\hhat_{\delta\gbar}^{(l)},\:
r^{(l-1)},\: l<k\:).$}
\label{eq:obshhat}
\end{equation}
On the other hand applying $\nabla_N^{k-1}$ to (\ref{eq:einsteinii2}) and 
evaluating on $M$ gives 
\begin{equation}
(k - 2)r^{(k-1)}|_M =
-\frac{k - 1}{n + 2}
h^{\delta\gbar}\hhat_{\delta\gbar}^{(k)}|_M\: + \: \textnormal{(terms involving
$\hhat_{\delta\gbar}^{(l)},\:r^{(l-1)},\: l<k\:).$}
\label{eq:obsr} 
\end{equation}
Substituting (\ref{eq:obsr}) into (\ref{eq:obshhat}) gives
\begin{equation}
\{k - (n+2)\}h^{\delta\gbar}
\hhat_{\delta\gbar}^{(k)}|_M = 
\: \textnormal{(terms involving
$\hhat_{\delta\gbar}^{(l)},\: l<k\:).$}
\label{eq:obs}
\end{equation}
Thus an obstruction to determining the trace of 
$\hhat_{\alpha\bbar}^{(k)}|_M$ occurs when the contents of the braces in 
(\ref{eq:obs}) is zero, i.e, when $k = n + 2$. Our inductive procedure will 
hence determine the trace of the Taylor coefficients of $\hhat_{\alpha\bbar}$ 
up to the coefficient of $\phi^{n+1}$.
On the other hand, applying $\nabla_N^k$ to (\ref{eq:einsteinab}) and 
evaluating on $M$ yields
$$
\hhat_{\alpha\bbar}^{(k)}|_M = -k\left(r^{(k-1)}
+ \frac{1}{n+2}h^{\delta\gbar}
\hhat_{\delta\gbar}^{(k)}\right)h_{\alpha\bbar}\Big|_M\: +
\: \textnormal{(terms involving
$\hhat_{\delta\gbar}^{(l)},\: l<k\:).$}
$$
Using (\ref{eq:obsr}) this shows that we can determine the trace-free part of 
the Taylor coefficients of $\hhat_{\alpha\bbar}$ up to the coefficient of 
$\phi^{n+2}$.

\subsection{Relation to Fefferman's approach} 
\label{subsec:fef}
Of course the fact that the obstruction is a scalar function
merely verifies Fefferman's result in~\cite{fef76}. In that paper it was shown
 that there is a smooth local solution $\rho$ to the Monge--Amp\'ere problem
\begin{equation}
\label{eq:fef}
\left\{ \begin{array}{ll} 
J(\rho)  := \det
\left(\begin{array}{cc}
\rho & \rho_{\overline{k}}\\
\rho_j & \rho_{j\overline{k}}\\
\end{array}\right)
 = -1 + O(\rho^{n+2}),\  \rho_j := \partial_j \rho, \textnormal{ etc}.; \\ 
\rho = 0\:\textnormal{on $M$}, \:\rho < 0\:\textnormal{on $X$}\\ 
\end{array} \right.
\end{equation}
that is unique modulo the addition of terms of order $O(\phi^{n+3})$.
For such a solution, Fefferman's approximately Einstein metric is given by 
the K\"ahler form 
$$
\omega = -\ddbar\log(-\rho(1 + O(\phi^{n+2}))),
$$
where the $O(\phi^{n+2})$ term is formally undetermined. Contracting with 
$W_\alpha \wedge W_{\bbar}$ and using (\ref{eq:ddbarphi}), 
$$
\omega(W_\alpha \wedge W_{\bbar}) = -\ddbar\log(-\rho)
(W_\alpha \wedge W_{\bbar}) + \phi^{n+1}vh_{\alpha\bbar} 
+ \textnormal{higher order terms in $\phi$},
$$
for a function $v$ on $M$ that is formally undetermined. In other words, an
ambiguity arises in the trace part of $\omega|_{\HH}$ at order $n+1$, just as
via our calculations.

\section{Volume renormalisation}
\label{sec:volume}

The Hermitian volume element $dv_{+}$ of $g_+$ on $X$ near $M$ is by definition
\begin{equation}
\label{eq:dvgw}
dv_{+} = \left(\frac{i}{2}\right)^{n+1}\frac{1}{(n+1)!}\omega^{n+1}. 
\end{equation}
Using (\ref{eq:diagmetric}), (\ref{eq:thetat}) and (\ref{eq:dtheta}) this is
\begin{equation}
\label{eq:dvg}
dv_{+} = \frac{(-1)^n}{2^{n+1}}\frac{\det\hhat}{n!\det h}\phi^{-n-2}\: 
d\phi\wedge\vartheta\wedge (d\vartheta)^n.
\end{equation}
We want to pull this form back to the product manifold 
$M\times (\delta,0]$ via our diffeomorphism $F$ (recall \S\ref{subsec:frame}). 
A basis for $T^\ast (M\times (\delta, 0])$ is 
$\{d\phi, \theta, \theta^\alpha,\theta^{\abar}\}$, where 
$\theta, \theta^\alpha,\theta^{\abar}$ are defined
 by pulling back $\vartheta, \vartheta^\alpha, 
\vartheta^{\abar}$ to $M$ under the inclusion $M\hookrightarrow\Ub$ and then 
extending 
constantly in the $N$-direction. For some uniquely determined 
function $s$ on $M\times (\delta,0]$ we have
$$
F^\ast(\vartheta\wedge (d\vartheta)^n)\equiv  s\: \theta\wedge (d\theta)^n
\mod d\phi.
$$
Contracting with the vector field $\partial/\partial\phi$, using  
(\ref{eq:dtheta}) and that  
$F_\ast(\partial/\partial\phi) = N$ and $\vartheta(N) = 0$ shows that in fact 
this is a strict equality:
$$
F^\ast(\vartheta\wedge (d\vartheta)^n) = s\: \theta\wedge (d\theta)^n.
$$
To determine $s$, take the exterior derivative $d$ of this equation and use
(\ref{eq:dtheta}) to give $s$ as the solution to the boundary-value problem
\begin{equation}
Ns = (n+1)rs\:;\quad\quad s|_M\equiv 1.
\label{eq:ivp}
\end{equation}
Since $\hhat$ is only known up to order $O(\phi^{n+1})$ we look for a formal 
solution to (\ref{eq:ivp}) up to order $O(\phi^{n+1})$. From (\ref{eq:dvg})
\begin{equation}
\label{eq:volpull5}
F^\ast(dv_{+}) = \frac{(-1)^n}{2^{n+1}}\frac{\det\hhat}{n!\det h}
s\phi^{-n-2}\: d\phi\wedge\theta\wedge (d\theta)^n.
\end{equation}
For some locally determined functions $v^{(j)}$ on $M$ this becomes
\begin{eqnarray}
\label{eq:volumeelement}
F^\ast(dv_{+}) &= &\phi^{-n-2}(v^{(0)} + v^{(1)}\phi + v^{(2)}\phi^2 +
\cdots +v^{(n+1)}\phi^{n+1} + \\
&&\quad\quad\quad\textnormal{\:\:higher order terms in $\phi$})\: d\phi
\wedge\theta\wedge (d\theta)^n.\nonumber
\end{eqnarray}

From here the procedure is identical to that in
the conformal setting (\cite{graham}). Pick a small number $\epsilon_0\in 
(\delta, 0)$, let $\epsilon_0<\epsilon<0$ and set $U^\epsilon$ to be the image 
under $F$ of $M\times (\epsilon_0, \epsilon)$. Write
\begin{eqnarray*} 
\Vol_{+}(\{\phi < \epsilon\}) &= &\const + \int_{U^\epsilon} dv_{+}\\
& = &\const + \int^\epsilon_{\epsilon_0}\!\int_M F^\ast(dv_{+}).
\end{eqnarray*}
Here ``$\const$'' is the volume of the compact set $\{\phi < \epsilon_0\}$ 
with 
respect to any Hermitian metric whose asymptotic expansion near $M$ agrees 
with the expansion of $g_+$.
It follows that 
\begin{equation}
\label{eq:volexpand}
\Vol_{+}(\{\phi < \epsilon\}) = c_0\epsilon^{-n-1} + c_1\epsilon^{-n}
+\cdots + c_n\epsilon^{-1} + L\log(-\epsilon) + V + o(1).
\end{equation}
The constant term $V$ is called the \textit{renormalised volume} in 
accordance with \cite{hs} and~\cite{graham}. 
The coefficients 
$c_j$ and $L$ are integrals over $M$ of local pseudohermitian invariants of 
$M$, with respect to the volume element $\theta\wedge (d\theta)^n$. In 
particular,
$$
L = \int_M v^{(n+1)}\theta\wedge (d\theta)^n.
$$
The following result is an analogue of Theorem 3.1 of~\cite{graham} (and is 
proved in the same way).
\begin{prop}
\label{prop:L}
The number $L$ is independent of the choice of pseudohermitian structure
 on $M$.
\end{prop}
\begin{proof}
Let $\theta$ and $\that = e^{2\Upsilon}\theta$ be two pseudohermitian 
structures on $M$, for $\Upsilon$ a function on $M$, with associated 
(by Lemma~\ref{lem:special}) special defining functions $\phi$ and $\phat$. So
$\phat = e^{2f(x,\phi)}\phi$, for a function $f$ in a neighbourhood of $M$, 
where we have used $x$ as a coordinate on $M$. 
We can inductively solve the equation $\phat = e^{2f(x,\phi)}\phi$ 
for $\phi$ to give $\phi = \phat b(x,\phat)$, for a uniquely determined 
positive function $b$. In this relation the $x$ still 
refers to the identification 
of $\Ub$ and $M\times [0,\delta)$ constructed using $\phi$. 
Set $\ehat(x,\epsilon) := \epsilon b(x,\epsilon)$. It follows that 
$\phat < \epsilon$ is 
equivalent to $\phi < \ehat(x,\epsilon)$, hence using (\ref{eq:volumeelement})
 we have
\begin{equation}
\begin{split}
\label{eq:voldiff}
\Vol_{+}(\{\phat < \epsilon\}) - \Vol_{+}(\{\phi < \epsilon\})
 &= \int_M \int_{\epsilon}^{\ehat}F^\ast(dv_{+})\\
 &\hspace{-55mm}= \int_M \int_{\epsilon}
^{\ehat}\phi^{-n-2}(v^{(0)} 
+ v^{(1)}\phi  + \cdots +v^{(n+1)}\phi^{n+1})\: d\phi\wedge\theta\wedge 
(d\theta)^n\ + o(1).
\end{split}
\end{equation}
In this expression the $\phi^{-1}$ term contributes
$\log b(x,\epsilon)$, so there is no $\log(-\epsilon)$ term as
$\epsilon\rightarrow 0$.
\end{proof}

While $L$ is independent of the choice of pseudohermitian structure, the 
renormalised volume is not. If $V_\theta$ and $V_{\that}$ are the renormalised
 volumes corresponding to pseudohermitian structures $\theta$ and $\that = 
e^{2\Upsilon}\theta$, then in accordance with~\cite{hs} and~\cite{graham} we
define the \textit{conformal anomaly} to be $V_{\that} - V_\theta$. So the
conformal anomaly is simply the constant term in the expansion 
(\ref{eq:voldiff}). From the form of the function $b$ used in the proof of 
Proposition~\ref{prop:L}, it follows that 
$$
V_{\that} - V_\theta = \int_M \PP_\theta(\Upsilon)\  \theta\wedge (d\theta)^n,
$$
where $\PP_\theta$ is a polynomial nonlinear differential operator whose
 coefficients are pseudohermitian invariants of $(M,\theta)$. It is 
important to note that the first variation of the conformal anomaly is 
\textit{not} simply $\int_M 2v^{(n+1)}\Upsilon\  \theta\wedge (d\theta)^n$, 
as one may initially suspect by analogy with the conformal setting 
of~\cite{graham}. The next section shows some of the extra computations and 
terms that arise.

\section{Renormalised volume when $n=1$ (proof of Theorem~\ref{thm:main})}
\label{sec:n=1}
To simplify the computation a little we may assume that 
$h^{1\oneb}\equiv 1$; this is
achieved through an appropriate choice of initial condition in the proof 
(equation (\ref{eq:ivp2})) of 
Lemma~\ref{lem:spframe}. We may now write $\hhat$ instead of $\hhat_{1\oneb}$ 
and $\Scal$ instead of $\Ric_{1\oneb}$.
Now applying $N^2$ to (\ref{eq:einsteinab}) we have
\begin{equation*}
\begin{split}
&3\hhat ^{\prime\prime} + 2\Big(3r^\prime  
+ \hhatinv^\prime(\hhat^\prime
- \frac{i}{2}\hhat_{,\:T})  +\hhatinv(\hhat^{\prime\prime}
- \frac{i}{2}\hhat_{,\:T N}) + \hhatinv_{,\:\overline{1}{ N}}\hhat_{,\:1}\\
&\quad\quad\quad\ \ + \hhatinv_{,\:\overline{1}}\hhat_{,\:1 N}
       + \hhatinv\hhat_{,\:1\overline{1} N}
       + \hhatinv^{\prime}\hhat_{,\:1\overline{1}} 
- \Scal^{\prime} \Big) = O(\phi).
\end{split}
\end{equation*}
Evaluating on $M$ and using the remark following Lemma~\ref{lem:comm3} gives
\begin{equation*}
\begin{split}
3\hhat^{\prime\prime} + 6r^\prime + 2\hhatinv^\prime\hhat^\prime  
+2\hhat^{\prime\prime}
- i\hhat^\prime_{,\:T}
+ 2\hhat^\prime_{,\:1\oneb}
       - 2\:\Scal^{\prime} = 0.
\end{split}
\end{equation*}
Using (\ref{eq:abbar}) this becomes
$$
2r^\prime =  -\frac{5}{3}\hhat^{\prime\prime} + \frac{2}{3}
(\hhat^\prime)^2 
+\frac{1}{3}\Delta_b \hhat^\prime
+ \frac{2}{3}\Scal^\prime,
$$
which after substituting (\ref{eq:hhatprime}) is
$$
2r^\prime = -\frac{5}{3}\hhat^{\prime\prime} 
+ \frac{1}{24}(\Scal)^2  
+\frac{1}{12}\Delta_b\Scal 
+ \frac{2}{3}\Scal^\prime.
$$
Finally using (\ref{eq:ricci}) yields
\begin{equation}
2r^\prime = -\frac{5}{3}\hhat^{\prime\prime}
+ \frac{1}{24}(\Scal)^2 - \frac{2}{3}|A|^2 
+\frac{1}{12}\Delta_b\Scal 
- \frac{2}{3}\Im A_{11,}^{\phantom{11,}11}\quad\textnormal{on $M$}.
\label{eq:rp}
\end{equation}
Turning now to equation (\ref{eq:einsteinii2}), a short computation, similar 
to the one just done, gives
\begin{equation}
\hhat^{\prime\prime} = 
\frac{1}{16}(\Scal)^2 - |A|^2 \quad\textnormal{on $M$}.
\label{eq:hhatpp}
\end{equation}
Therefore up to terms of order $O(\phi^2)$ we have that
\begin{equation}
\label{eq:hhatn=1}
\hhat = 1 + \frac{1}{4}\Scal\phi 
+ \frac{1}{2}\left(\frac{1}{16}(\Scal)^2 - |A|^2\right)\phi^2,
\end{equation}
where $\Scal$ and $A_{11}$ are evaluated on $M$.

When $n = 1$, it is easily seen that the solution up to the
second jet of the boundary-value problem (\ref{eq:ivp}) is
$$
s = 1 + (r^\prime|_M)\phi^2 + O(\phi^3).
$$
Then from (\ref{eq:volpull5}) and (\ref{eq:volumeelement}) 
the volume form on the product manifold $M \times (\delta, 0]$ is
\begin{equation}
\label{eq:voleln=1}
-\frac{1}{4}\phi^{-3}\Bigl(1 
+ (\hhat^{\prime}|_M)\phi + 
\frac{1}{2}(\hhat^{\prime\prime}|_M)\phi^2 + \cdots\Bigr)
\Bigl(1 + (r^\prime|_M)\phi^2 + \cdots\Bigr)\ d\phi\wedge\theta\wedge d\theta.
\end{equation}
The coefficient of $\phi^{-1}$ in this expansion is 
\begin{equation*}
\begin{split}
v^{(2)} &= -\frac{1}{8}(2r^\prime 
+ \hhat^{\prime\prime})|_M\\
        &= -\frac{1}{8}\left(\frac{1}{12}\Delta_b\Scal 
- \frac{2}{3}\Im A_{11,}^{\phantom{11,}11} \right)\Big|_M,
\end{split}
\end{equation*}
using (\ref{eq:rp}) and (\ref{eq:hhatpp}). Therefore, 
\begin{equation*}
\begin{split}
L &= -\frac{1}{8}\int_M\frac{1}{12}\Delta_b\Scal 
- \frac{2}{3}\Im A_{11,}^{\phantom{11,}11} \ \theta\wedge d\theta\\
  & = 0,
\end{split}
\end{equation*}
by the divergence formula (\ref{eq:div}).

We now proceed to calculate the conformal anomaly, using the notation in the 
proof of Proposition~\ref{prop:L}. By (\ref{eq:voldiff}) the conformal anomaly
is the constant term in $\epsilon$ in
$$
\int_M\int_{\epsilon}
^{\ehat}\phi^{-3}(v^{(0)} + v^{(1)}\phi +v^{(2)}\phi^2)
\: d\phi\wedge\theta\wedge d\theta.
$$
Performing the integration with respect to $\phi$ this constant term 
eventually simplifies to 
$$
\int_M -v^{(0)}f^{\prime\prime}|_M  -2v^{(1)}f^{\prime}|_M -2v^{(2)}\Upsilon
\:\theta\wedge d\theta,
$$
where we recall that $\phih = e^{2f}\phi$ is the special defining function 
associated to $\that$. Using (\ref{eq:volumeelement}) and (\ref{eq:voleln=1}) 
we calculate the
coefficients $v^{(j)}$ to get the conformal anomaly as
\begin{equation}
\label{eq:ca2}
V_{\that} - V_\theta = \int_M \frac{1}{4}f^{\prime\prime}|_M  
+\frac{1}{8}\Scal f^{\prime}|_M 
+ \frac{1}{48}(\Delta_b\Scal - 8\Im A_{11,}^{\phantom{11,}11})\Upsilon
\:\theta\wedge d\theta.
\end{equation}

Now from the proof of Lemma~\ref{lem:special}, the function $f$ is the unique
solution to the boundary-value problem given by the PDE
$$
\frac{4}{\phi}\Re \langle\partial\phi,\partial f\rangle_{g_+} + 
4|\partial f|^2_{g_+} = 0
$$ 
and the boundary condition $f|_M = \Upsilon$. The PDE becomes
\begin{equation}
\label{eq:f}
\begin{split}
&2f^{\prime} - \phi(f_1\eta^1 + f_{\oneb}\eta^{\oneb}) - 2\hhatinv^{1\oneb}
f_{1}f_{\oneb} + 2\phi^3\eta^1\eta^{\oneb}f_1 f_{\oneb} + 2\phi f_0f_{\zerob}\\
&\quad -2\phi^2\eta^1f_1f_{\zerob} - 2\phi^2\eta^{\oneb}f_{\oneb}f_0 = 0.
\end{split}
\end{equation}
Thus $f^{\prime}|_M = \Upsilon_1\Upsilon^1$. 

To compute $f^{\prime\prime}|_M$ we,
as usual, apply $\nabla_N$ to (\ref{eq:f}) while working in a special frame
given by Lemma~\ref{lem:spframe} and then evaluate on $M$. We also need to 
commute covariant derivatives of the function $f$, for which we use the 
formula, derived similarly to those in Lemma~\ref{lem:comm2},
$$
f_{1N} = f_{N1} - \frac{1}{2}f_1 r 
- \frac{i}{2}f_{\oneb}A^{\oneb}_{\phantom{\oneb}1}.
$$ 
We obtain
\begin{equation*}
\begin{split}
f^{\prime\prime}|_M = & \frac{1}{2}(\Upsilon_1\eta^1 
+ \Upsilon_{\oneb}\eta^{\oneb}) - (\Upsilon_1\Upsilon^1)^2 
- \frac{1}{4}(\Upsilon_T)^2 - \frac{1}{4}\Scal\Upsilon_1\Upsilon^1\\ 
&+ 2\Re(\Upsilon_{11}\Upsilon^1\Upsilon^1) 
- \Upsilon_1\Upsilon^1\Delta_b\Upsilon - \Im(\Upsilon_1\Upsilon_1 A^{11}).
\end{split}
\end{equation*}
From (\ref{eq:einsteinai}), 
\begin{equation*}
\eta^1 = \frac{1}{3}\left(-\frac{1}{4}\Scal^1 
- iA_{\phantom{1}\oneb,}^{1\phantom{1\oneb,}\oneb}\right).
\end{equation*}
Substituting the previous two formulae into (\ref{eq:ca2}) and integrating by
parts gives the formula for the conformal anomaly in Theorem~\ref{thm:main}, 
thereby completing the proof.

\section{A renormalised Chern--Gauss--Bonnet formula when $n=1$ (proof 
of Corollary~\ref{cor:main})}
\label{sec:cgb}

Recall from the statement of Corollary~\ref{cor:main}
that we now suppose we are given a complete Einstein--K\"ahler metric 
$g_{\textnormal{\tiny{EK}}}$ on $X$. Let $X^\epsilon$ denote the region 
enclosed by 
the hypersurface $M^\epsilon$.
The Chern--Gauss--Bonnet formula (\cite{chern}) for a complex 
manifold-with-boundary tells us that
\begin{equation}
\label{eq:cgbk}
\chi(X^\epsilon) = \int_{X^\epsilon} c_2 + \int_{M^\epsilon} \Pi^\epsilon,
\end{equation}
where $c_2$ is the second Chern form and $\Pi^\epsilon$ is the restriction
to $M^\epsilon$ of a certain top-form $\Pi$ on $X$. The original formula
 for $\Pi$  consists of wedge products of Levi-Civita curvature and 
connection components in an orthonormal basis. It is straightforward to derive
 from this an expression in basis-free notation. Indeed, let $\nu$ be the unit
outward normal to $M^\epsilon$ and $R$ denote the curvature of 
$X$, viewed as a two-vector-valued two-form. If $\Sigma = Z\otimes\zeta$ is a 
$k$-vector-valued $k$-form and $\Sigma^\prime = Z^\prime\otimes\zeta^\prime$ 
is an $k^\prime$-vector-valued $k^\prime$-form then 
$\Sigma\wedge \Sigma^\prime$ will denote the 
$(k+k^\prime)$-vector-valued $(k+k^\prime)$-form 
$Z\wedge Z^\prime\otimes\zeta\wedge
\zeta^\prime$.
Then it holds that
\begin{equation}
\label{eq:pie}
\Pi^\epsilon = \frac{1}{12\pi^2}dv_{\textnormal{\tiny{EK}}}^{\epsilon}\hook
(D\nu\wedge D\nu\wedge D\nu
+ 3 D\nu\wedge R),
\end{equation}
where for $\Sigma$ a three-vector-valued three-form on $M^\epsilon$, 
$dv_{\textnormal{\tiny{EK}}}^{\epsilon}\hook(\Sigma)$ denotes the 
three-form given by 
contraction of $\Sigma$ with the volume form on $M^\epsilon$ 
induced from that on $(X, g_{\textnormal{\tiny{EK}}})$. (In the above, we are
assuming skew-symmetry as vectors as well as forms.)

Now as $\epsilon\to 0$, the integral over $X^\epsilon$ in (\ref{eq:cgbk}) 
will diverge; however 
Burns--Epstein~\cite{be} proved that replacing $c_2$ by the 
renormalised Chern form $\tilde{c_2} := c_2 - (1/3)c_1^{\phantom{1}2}$ 
results in a convergent integral. Thus we may rewrite (\ref{eq:cgbk}) as
\begin{equation}
\label{eq:cgbk2}
\chi(X^\epsilon) = \int_{X^\epsilon}\left(c_2 - \frac{1}{3}c_1^{\phantom{1}2}
\right) +
\frac{1}{3}\int_{X^\epsilon}c_1^{\phantom{1}2} + \int_{M^\epsilon}\Pi^\epsilon,
\end{equation}
where the first integral is now convergent. If we expand the second and third 
integrals as power series in $\epsilon$, it follows that as 
$\epsilon\to 0$ only the 
constant term in these series will remain. By the Einstein condition,
\begin{equation*}
c_1^{\phantom{1}2} = \left(-(n+2)\frac{i}{2\pi}\omega_{\textnormal{\tiny{EK}}}
\right)^2 = \frac{18}{\pi^2}dv_{\textnormal{\tiny{EK}}},
\end{equation*}
where $\omega_{\textnormal{\tiny{EK}}}$ ($dv_{\textnormal{\tiny{EK}}}$) is 
the K\"ahler (volume) form of 
$g_{\textnormal{\tiny{EK}}}$. Letting $\epsilon\to 0$ in (\ref{eq:cgbk2}) 
we obtain
\begin{equation}
\label{eq:cgbk3}
\chi(X) = \int_{X}\left(c_2 - \frac{1}{3}c_1^{\phantom{1}2}\right) +
\frac{6}{\pi^2}V_\theta + \int_M S_\theta \:\theta\wedge d\theta,
\end{equation}
where $V_\theta$ is the renormalised volume with respect to 
the pseudohermitian structure $\theta$ and $S_\theta\ \theta\wedge d\theta$ 
is defined to be the pullback to $M$ of the constant term in the power
 series expansion of $\Pi^\epsilon$.

In order to compute $S_\theta$ one could use the explicit expression for 
$\Pi^\epsilon$ and consider its asymptotics as $\epsilon\to 0$. This 
computationally tedious 
approach though can be avoided by using some elementary invariant theory and 
the 
formula obtained earlier for the conformal anomaly. We first observe how 
$S_\theta$ changes under a particular change of pseudohermitian structure. 
\begin{lem}
\label{lem:Ss}
If $\that = e^{2c}\theta$, for a 
constant $c$, then $S_{\that} = e^{-4c}S_{\theta}$.
\end{lem}
\begin{proof}
From Lemma~\ref{lem:special}, the special defining functions $\phi$ and $\phat$
associated to the pseudohermitian structures $\theta$ and $\that$, 
respectively, satisfy the relation $\phat = e^{2c}\phi$. Hence $\phi$ and 
$\phih$ have the same level sets, implying
$$
S_\theta \: \theta \wedge d\theta = S_{\that} \: \that \wedge d\that
$$
and completing the proof.
\end{proof}

We can now see the form that $S_\theta$ must take:
\begin{lem}
\label{lem:ab}
There exist universal constants $a,b$ such that
$$
S_\theta = a(\Scal)^2 + b|A|^2 + \textnormal{ divergence}.
$$
Here ``divergence'' denotes a function whose integral over $M$ will vanish via 
the divergence formula (\ref{eq:div}).
\end{lem}
\begin{proof}
By (\ref{eq:pie}) it is clear that $S_\theta$ must be a polynomial expression
in the components of the pseudohermitian curvature and torsion and their 
pseudohermitian covariant derivatives. It also must hold that $S_\theta$ be
 invariant under a change of frame for the holomorphic tangent bundle. In 
other words, in the terminology of~\cite{hirachi}, $S_\theta$ is a scalar
pseudohermitian invariant. But the scalar pseudohermitian invariants that
satisfy the same transformation law as $S_\theta$ does in Lemma~\ref{lem:Ss} 
are shown in the proof of Theorem 5.1 of~\cite{hirachi} to be
$$
\Scal_1^{\phantom{1}1},\  \Scal_{\oneb}^{\phantom{\oneb}\oneb},\ 
A_{11,}^{\phantom{11,}11},\  
A_{\oneb\oneb,}^{\phantom{\oneb\oneb,}\oneb\oneb},\ 
(\Scal)^2, \ |A|^2,
$$
and the result follows.
\end{proof}

To find the constants in Lemma~\ref{lem:ab}, we compute the first variation 
of the conformal anomaly using (\ref{eq:cgbk3}) and compare it with the value
given by Theorem~\ref{thm:main}. It suffices to do this for a special kind of 
pseudohermitian manifold, since the constants in question are universal. 
We will need the following basic result; for a proof see, e.g.,~\cite{lee86}.
\begin{lem}
\label{lem:lee} 
If $\widehat{A}_{\alpha\beta}$ and $\widehat{\Scal}$ denote the
pseudohermitian torsion and scalar curvature of $M$ with respect to the
pseudohermitian structure $\that = e^{2\Upsilon}\theta$, for a function 
$\Upsilon$, then
$$
\hat{A}_{\alpha\beta} = e^{-2\Upsilon}(A_{\alpha\beta} + 2i\Upsilon
_{\alpha\beta} - 4i\Upsilon_\alpha\Upsilon_\beta),
$$
$$
\widehat{\Scal} = e^{-2\Upsilon}(\Scal + 4\Delta_b\Upsilon 
- 8\Upsilon_\alpha\Upsilon^\alpha).
$$
\end{lem}
Now in general, from (\ref{eq:cgbk3}), the quantity
$$
\VV(M) := \frac{6}{\pi^2}V_\theta + \int_M S_\theta \:\theta\wedge d\theta
$$
is a global CR invariant. Let $\hat{\Scal}_t, \hat{A}_t$ and $V_{\that_t}$ 
denote the scalar curvature, torsion and renormalised
 volume with respect to the pseudohermitian structure $\that_t := 
e^{2\Upsilon t}\theta$, for a real parameter $t$. Then
\begin{equation}
\begin{split}
\label{eq:ca}
V_{\that_t} - V_\theta &= \frac{\pi^2}{6}\left(
\int_M S_\theta \:\theta\wedge d\theta - \int_M S_{\that_t} \:\that_t\wedge 
d\that_t\right)\\
& = \frac{\pi^2}{6}
\int_M a\Big((\Scal)^2 - (e^{2t\Upsilon}\widehat{\Scal_t})^2\Big) 
+ b\Big(|A|^2 - |e^{2t\Upsilon}\hat{A}_t|^2\Big)\  \theta\wedge d\theta\\
& = -\frac{\pi^2}{3}
\int_M t(4a\Upsilon\Delta_b\Scal + 2b\Upsilon\Im A_{11,}^{\phantom{11,}11})
\ \theta\wedge d\theta + \int_M O(t^2);
\end{split}
\end{equation}
the last equality here is obtained from Lemma~\ref{lem:lee} and integrating 
by parts. The first variation of the conformal anomaly is 
$$
\frac{d}{dt}V_{\that_t}|_{t=0} = -\frac{\pi^2}{3}
\int_M 4a\Upsilon\Delta_b\Scal + 2b\Upsilon\Im A_{11,}^{\phantom{11,}11}
\ \theta\wedge d\theta.
$$

It is a fact that there exists a CR three-manifold, 
embeddable in a complex two-manifold, 
with a pseudohermitian structure whose torsion is identically zero but whose
scalar curvature is not everywhere-annihilated by the sublaplacian. Let now
$(M,\theta)$ be this manifold.
Set $\Upsilon = \Delta_b\Scal$.  Then 
$$
\frac{d}{dt}V_{\that_t}|_{t=0} = -\frac{4a\pi^2}{3}
\int_M (\Delta_b\Scal)^2 \ \theta\wedge d\theta.
$$
But in a similar way, Theorem~\ref{thm:main} tells us that
$$
\frac{d}{dt}V_{\that_t}|_{t=0} = \frac{1}{96}
\int_M (\Delta_b\Scal)^2 \ \theta\wedge d\theta.
$$
We conclude that 
$$
a = -\frac{1}{128\pi^2}.
$$

Returning to the general case, this time we set 
$\Upsilon = \Im A_{11,}^{\phantom{11,}11}$---it is a fact 
that there exist pseudohermitian manifolds on which this function is 
not identically zero---and proceed as above to obtain
$$
\frac{d}{dt}V_{\that_t}|_{t=0} = \int_M \frac{1}{96}
\Im A_{11,}^{\phantom{11,}11}\Delta_b\Scal 
- \frac{2b\pi^2}{3}(
\Im A_{11,}^{\phantom{11,}11})^2 \ \theta\wedge d\theta.
$$
But this time Theorem~\ref{thm:main} tells us that
$$
\frac{d}{dt}V_{\that_t}|_{t=0} = \int_M \frac{1}{96}
\Im A_{11,}^{\phantom{11,}11}\Delta_b\Scal 
- \frac{1}{12}(
\Im A_{11,}^{\phantom{11,}11})^2 \ \theta\wedge d\theta.
$$
We conclude that 
$$
b = \frac{1}{8\pi^2},
$$
and Corollary~\ref{cor:main} is proved.

\begin{rem}
A similar Chern--Gauss--Bonnet formula will hold if $X$ is any complex manifold
with Hermitian metric $g$ whose asymptotic expansion near $M$ agrees with 
that of Fefferman's approximately Einstein metric $g_+$. For 
simplicity, if we choose the Hermitian metric used in the definition of 
renormalised volume (\S~\ref{sec:volume}) to be $g$, then for a small 
negative number $\epsilon_0$, if $\Vol_g(\{\phi < \epsilon_0\})$ denotes 
the volume of the set $\{\phi < \epsilon_0\}$ with respect to $g$, we have
$$
\chi(X) = \int_X \left(c_2 - \frac{1}{3}c_1^{\phantom{1}2}\right) + \VV(M) 
- \frac{6}
{\pi^2}
\Vol_g(\{\phi < \epsilon_0\}) + \frac{1}{3}\int_{\{\phi<\epsilon_0\}}c_1^2\:.
$$
\end{rem}

\section{Example: The Bergman metric on the unit ball in $\C^2$}
\label{sec:bergman}

Let $\Bbar$ denote the closed unit ball in $\C^2\ni(z,w)$. Its interior $B$ 
carries a unique globally defined Einstein--K\"ahler metric, the Bergman 
metric, with K\"ahler
form 
\begin{equation}
\label{eq:Bergman}
\omega_{\textnormal{\tiny{B}}} = \ddbar\log\left(-\frac{1}{\rho}\right),
\end{equation}
where
$$
\rho(z,w) = R^2 - 1
$$
and
$$
R := \sqrt{z\zbar + w\wbar}.
$$

The unit sphere $S$ forming the boundary of $\Bbar$ inherits a natural CR
structure from $\C^2$. The so-called standard pseudohermitian structure on 
$S$ is
$$
\theta = \frac{i}{2}(\dbar\rho - \partial\rho)|_S = -i\partial\rho|_S
= -i(\zbar dz + \wbar dw)|_S.
$$
For this $\theta$,
$$
\Scal \equiv 2;\quad A \equiv 0
$$
(\cite{webster}). It is also easy to check that
$$
\theta\wedge d\theta = 2 dv_{\textnormal{\tiny{S}}},
$$
where $dv_{\textnormal{\tiny{S}}}$ denotes the usual volume form on $S$ given 
by restriction of the Euclidean volume form on $\C^2$.

Our goals in this section are to compute the renormalised volume $V$ of 
$B$ with respect to $\omega$ and $\theta$ and to verify the
renormalised Chern--Gauss--Bonnet formula in Corollary~\ref{cor:main}.

The first step is to find the special defining function $\phi$ associated 
to $\theta$. In fact one may verify that
\begin{equation}
\label{eq:specialphi}
\phi = 4\left(\frac{R - 1}{R + 1}\right)
\end{equation}
solves the boundary-value problem in Lemma~\ref{lem:special} for 
$\theta$ as above. (Suspecting a defining function of this form was motivated 
by a similar defining function appearing in the conformal setting for the 
Poincar\'e metric on the ball (\cite{graham})).

Now in general the asymptotic procedures of the preceeding sections do not
allow one to compute the renormalised volume. However in our particular
example we will be able to compute $V$ directly by ad hoc methods. Indeed, 
a straightforward calculation gives the 
volume form for the Bergman metric as
$$
dv_{\textnormal{\tiny{B}}} = 
\frac{1}{(1-R^2)^3}dv_{\textnormal{\tiny{Eu}}},
$$
where $dv_{\textnormal{\tiny{Eu}}}$ denotes the Euclidean volume form. In 
polar coordinates
$$
dv_{\textnormal{\tiny{B}}} = 
\frac{R^3}{(1-R^2)^3}dR\wedge dv_{\textnormal{\tiny{S}}}.
$$
Finally in terms of the special coordinate $\phi$ we obtain
\begin{equation}
\label{eq:volformB}
dv_{\textnormal{\tiny{B}}} = 
\phi^{-3}\left(-\frac{1}{2} - \frac{1}{4}\phi + \frac{1}{64}\phi^3 + 
\frac{1}{512}\phi^4\right)d\phi\wedge dv_{\textnormal{\tiny{S}}}.
\end{equation}
Thus 
$$
\Vol_{\textnormal{\tiny{B}}}(\{\phi < \epsilon\}) = 
\int_S dv_{\textnormal{\tiny{S}}}
\int^\epsilon_{-4} \:\phi^{-3}\left(-\frac{1}{2} - \frac{1}{4}\phi 
+ \frac{1}{64}\phi^3 + \frac{1}{512}\phi^4\right)d\phi.
$$
Integrating and looking at the constant term in $\epsilon$ we get
$$
V = \frac{3\pi^2}{16}.
$$
Observe also that there is no log term in $\epsilon$, as we should expect.

We now present a second derivation of the formula (\ref{eq:volformB}) that 
is longer but will illustrate concretely most of the rather technical frame 
change calculations and the like of the preceeding sections.
Recall from the formula (\ref{eq:volpull5}) for the volume form and the 
subsequent discussion that the ingredients we need are \textit{exact} formulae
(i.e., not just asymptotic approximations---we are interested in the 
constant term $V$) for the following: the tensor $\hhat$ appearing 
in the diagonal expression (\ref{eq:diagmetric}) for the K\"ahler form; the 
tensor $h$
appearing in (\ref{eq:ddbarphi}); and the function $s$ that solves 
(\ref{eq:ivp}). 
We shall moreover work with a special frame given by Lemma~\ref{lem:spframe} 
that in addition satisfies $h_{1\oneb}|_S\equiv 1$. By working with such a 
frame we may ignore the $h$ term appearing in (\ref{eq:volpull5}) and also 
more 
readily verify the formulae given for $\hhat$ and the transverse curvature 
$r$ in \S~\ref{sec:n=1}. In other words we shall, for our particular example, 
carry out all of 
the computations in the volume renormalisation procedure, with the exception 
of solving the Einstein equation for $\hhat$. The formula for $\hhat$ will 
instead be directly apparent.

Let us begin with the appropriate choice of frame. The natural CR structure on
$S$ is usually defined by the restriction to $S$ of 
$$
\vartheta^1 = -wdz + zdw.
$$
For this coframe in fact $h_{1\oneb}|_S\equiv 1$. The dual frame is 
$$
W_1 = \frac{1}{R^2}\left(-\wbar\frac{\partial}{\partial z} 
+ \zbar\frac{\partial}{\partial w}\right).
$$
Then with special $\phi$ defined as in (\ref{eq:specialphi}), we compute
the vector field $\xi$ from its characterisation (\ref{eq:xi}):
$$
\xi = \frac{(R+1)^2}{4R}\left(z\frac{\partial}{\partial z} + w\frac{\partial}
{\partial w}  \right).
$$

Now if we take $\varthetat^1 = \vartheta^1$ then it is readily observed that
$$
\hhat_{1\oneb} = \frac{4}{R^2(R+1)^2}
$$
satisfies (\ref{eq:diagmetric}). However we want to work with a special frame
given by Lemma~\ref{lem:spframe}; thus we need to multiply $\varthetat^1$ by 
the function $a$ that solves (\ref{eq:ivp2}). The connection form 
$\phi_{1}^{\phantom{1}1}$ may be determined from the structure
equation (\ref{eq:glstr}). We obtain
$$
\phi_{1\phantom{1}N}^{\phantom{1}1} = \frac{-5R^2 - 8R -3}{16R}
$$
and then
$$
a = \frac{2}{R^{3/2}(R+1)}.
$$
Therefore in this special frame
$$
\hhat = \dfrac{\left(\dfrac{4}{R^2(R+1)^2}\right)}
{\left(\dfrac{4}{R^3(R+1)^2}\right)} = R.
$$
In terms of $\phi$,
$$
\hhat = \frac{4 + \phi}{4 - \phi}.
$$
Expanding this as a Taylor series in $\phi$, we observe that it agrees with 
(\ref{eq:hhatn=1}) up to the coefficient of $\phi^2$.

The transverse curvature, which was defined as $r = 2\ddbar\phi(\xi\wedge
\overline{\xi})$, is
$$
r = \frac{1 - R^2}{8R}.
$$
Note that $r$ vanishes on the boundary. In terms of $\phi$,
$$
r = \frac{2\phi}{\phi^2 -16}.
$$
Observe that $r^{\prime}|_S = -(1/8)$, agreeing with (\ref{eq:rp}). We now get
that
$$
s = \frac{1}{256}(\phi^2 - 16)^2
$$
is the solution to the boundary-value problem (\ref{eq:ivp}).

Substituting the above formulae for $\hhat$ and $s$ into the formula 
(\ref{eq:volpull5}) for the volume form we get 
(\ref{eq:volformB}).

Let us now verify the renormalised Chern--Gauss--Bonnet formula in 
Corollary~\ref{cor:main}. Since $\Scal\equiv 2, A\equiv 0$ and $V = 3\pi^2/16$,
we have that the CR invariant $\VV(S) = 1$. The Levi-Civita 
curvature matrix and hence Chern forms of the Bergman metric may be 
computed from (\ref{eq:Bergman}). It turns out that 
$c_2 - (1/3)c_1^2\equiv 0$. Since the Euler characteristic of the ball 
$\chi(B) = 1$, the verification is complete. Remark that this is moreover 
consistent with the value $\mu(S) = -1$ given for the Burns--Epstein
invariant in~\cite{be88}.

\appendix
\section{CR $Q$-curvature}
\label{sec:Q}

A straightforward adaptation of the method in~\cite{fg} to our CR setting will
 lead to the proof of (\ref{eq:CRQ}). We begin by recalling a result from that 
paper. To 
distinguish notation between the conformal and CR cases, bold characters are 
used in the former. 
\begin{thm}[{\cite[Theorem 3.1]{fg}}]
\label{thm:fg}
Let $(\Mb,[\gb])$ be a smooth compact conformal manifold of 
even dimension $\nb$.  Let $\gbp$ 
be the Poincar\'e metric on $\Xbo := \Mb\times (-1,0]$. Choose a 
representative metric 
$\gb$ with special defining function $\pb$ that satisfies 
$\pb^2\gbp|_{T\Mb} = \gb$, $|d(\log(-\pb))|_{\gbp}^2 = 1$. 
There is a unique solution $\Ubo$ 
mod $O(\pb^{\nb})$ to 
$$
\Db_{\gbp}\Ubo = \nb + O(\pb^{\nb+1}\log(-\pb))
$$ 
of the form 
$$
\Ubo = \log(-\pb) + \Ab + \Bb\pb^{\nb}\log(-\pb) + O(\pb^{\nb})
$$
with 
$$
\Ab,\Bb\in C^ \infty(\Xbo),\quad \Ab|_{\Mb} = 0.
$$
Also, $\Ab$ mod $O(\pb^{\nb})$ and $\Bb|_{\Mb}$ are formally determined by 
$\gb$.
\end{thm}
The $Q$-curvature is then defined to be a constant 
multiple of of $\Bb|_{\Mb}$. It is further shown in~\cite{fg} that the 
integral of this quantity is a constant multiple of the log term coefficient
in the volume expansion of $(\Xbo, \gbp)$. Returning to CR geometry, we have
the following analogue of Theorem~\ref{thm:fg}. 
\begin{prop}
\label{prop:lap}
Let $M$ be a smooth compact strictly pseudoconvex CR manifold of 
dimension $2n+1$ that forms the boundary of a complex manifold $X$. Let $g_+$ 
be Fefferman's approximately Einstein metric on $X$ near $M$. 
Choose any 
pseudohermitian structure $\theta$ for $M$ and let $\phi$ be the special 
defining function associated to $\theta$ according to Lemma~\ref{lem:special}.
There is a unique solution $U$ mod $O(\phi^{n+1})$ to 
$$
\Delta_{g_+}U = 2(n+1) + O(\phi^{n+2}\log(-\phi))
$$ 
of the form 
$$
U = \log(-\phi) + A + B\phi^{n+1}\log(-\phi) + O(\phi^{n+1})
$$
with 
$$
A,B\in C^ \infty(X),\quad A|_M = 0.
$$
Also, $A$ mod $O(\phi^{n+1})$ and $B|_M$ are formally determined by $\theta$.
\end{prop}
\begin{proof}
In~\cite{gl} the  K\"ahler Laplacian 
$\Delta_{g_+}$ was decomposed into tangential and normal pieces 
with respect to the foliation $\{\rho = \epsilon\}$, where $\rho$ is a
local solution to the Monge--Amp\'ere problem (\ref{eq:fef}). We, however,
require a decomposition with respect to the foliation $\{M^\epsilon\}
= \{\phi = \epsilon\}$.
To this end, we note that from the form (\ref{eq:diagmetric}) of $g_+$ 
we may write, for a function $f$ on $X$,
\begin{equation*}
\begin{split}
\Delta_{g_+}f&= - \left(\ginv^{ 0\zerob}(D_{\xit} D_{\overline{\xit}}
                                   +D_{\overline{\xit}} D_{\xit})
                + \ginv^{\alpha\bbar}(D_{W_\alpha} D_{W_{\bbar}}
                                     +D_{W_{\bbar}} D_{W_\alpha})\right) f\\
             &= -\left(\phi^2(D_{\xit} D_{\overline{\xit}}
                              +D_{\overline{\xit}} D_{\xit})-
                \phi\hhatinv^{\alpha\bbar}(D_{W_\alpha} D_{W_{\bbar}}
                                       +D_{W_{\bbar}} D_{W_\alpha})\right) f .
\end{split}
\end{equation*}
Now by the definition of covariant differentiation
\begin{equation*}
\begin{split}
D_{\xit} D_{\overline{\xit}} f &= \xit\;\ \!\!\xitbar f 
- f_{\gbar}\overline{\psi_{ 0\phantom{\gamma}\zerob}^
{\phantom{ 0}\gamma}} - (\xitbar f)\overline{\psi_{ 0\phantom{\gamma}\zerob}
^{\phantom{ 0} 0}}\\
                  &= N^2 f   
+ (r + \phi^2 |\eta|^2)N f + Y_1 f,
\end{split}
\end{equation*}
similarly,
\begin{equation*}
D_{\overline{\xit}} D_{\xit}f = N^2 f + (r + \phi^2 |\eta|^2)N f + Y_2 f,
\end{equation*}
and
\begin{equation*}
D_{W_\alpha} D_{W_{\bbar}}f = (Nf)h_{\alpha\bbar} + Y_3f,
\end{equation*}
\begin{equation*}
D_{W_{\bbar}}D_{W_\alpha} f = (Nf)h_{\alpha\bbar} + Y_4f,
\end{equation*}
for tangential operators $Y_1, Y_2, Y_3, Y_4$.
Collecting these formulae together we have that
\begin{equation}
\label{eq:laplacian}
\Delta_{g_+} = -2\phi^2N^2
- 2\phi(\phi r + \phi^3 |\eta|^2 -  \hhatinv_\alpha^{\phantom{\alpha}\alpha})
N  +\phi Y,
\end{equation}
where $Y$ is a tangential operator.

The above formula implies that if $a_j \in C^ \infty(M)$ then 
\begin{equation*}
\begin{split}
\frac{1}{2}\lap(a_j\phi^j) &= -j(j-1)a_j\phi^j 
+ j\hhatinv_\alpha^{\phantom{\alpha}\alpha}a_j\phi^j 
+ O(\phi^{j+1})\\  
                &= -j(j-(n+1))a_j\phi^j + O(\phi^{j+1}).
\end{split}
\end{equation*}
Hence by setting $a_0 = 0$ we can inductively determine $a_1,\dots,a_n$ 
uniquely so that, taking $A = \sum_{j = 0}^n a_j\phi^j$, the function 
$\log(-\phi) + A$ solves 
$$
\frac{1}{2}\Delta_{g_+}(\log(-\phi) + A) = n + 1 + E\phi^{n+1}.
$$
Here $E\in C^\infty(X)$ and $E|_M$ is formally determined by the 
pseudohermitian structure $\theta$. From (\ref{eq:laplacian}) it is easy to 
compute that if $B\in C^\infty(X)$,
\begin{equation}
\label{eq:B}
\frac{1}{2}\lap(B\phi^{n+1}\log(-\phi)) = -(n+1)B\phi^{n+1} + O(\phi^{n+2}
\log(-\phi)).
\end{equation}
So setting $B = \frac{1}{n+1}E$ completes the proof of the proposition, since 
the uniqueness assertion is clear from construction.
\end{proof}

\begin{prop} If $L$ denotes the log term coefficient in the asymptotic 
expansion of the volume of $(X,g_+)$ and $B$ is as in Proposition
~\ref{prop:lap} then
$$
L = \frac{(-1)^{n+1}}{2^{n+1}n!}\int_M B \ \theta\wedge (d\theta)^n.
$$
\end{prop}
\begin{proof} One can check that the unit outward normal to each $M^\epsilon$
is $\sqrt{2}\phi\Re\:\xit$. Taking the function $U$ from 
Proposition~\ref{prop:lap}, 
for a small negative number $\epsilon_0$ Green's theorem gives
$$
\int_{\epsilon_0 <\phi < \epsilon }\lap U dv_{+}
 = -\sqrt{2}\epsilon\int_{M^\epsilon} (\Re\:\xit) U + \sqrt{2}\epsilon_0\int_
{M^{\epsilon_0}}(\Re\:\xit) U,
$$
where the integrals on the right-hand side are with respect to the volume 
element induced on the each respective hypersurface. We want to
compare coefficients of $\log(-\epsilon)$ in this equation. By 
Proposition~\ref{prop:lap} and (\ref{eq:volexpand})
 the coefficient of $\log(-\epsilon)$ on the left-hand side is $2(n+1)L$. 
As for the
right-hand side, only the first integral is relevant. We use our
diffeomorphism $F$ to pull this integral back to $M$. Together with the form
of $U$ given in Proposition~\ref{prop:lap} and the expression for the volume 
element (\ref{eq:volumeelement}), this integral becomes
\begin{equation*}
\begin{split}
&-\sqrt{2}\epsilon\int_M \Big(\epsilon^{-1} + (\Re\:\xit) A 
+ (n + 1)B\epsilon^{n}\log(-\epsilon) + O(\epsilon^n)\Big)\times\\
&\quad\quad\quad\quad \sqrt{2}\epsilon^{-n-1}\Big(v^{(0)} + v^{(1)}\epsilon
+\cdots
\Big) 
\ \theta\wedge (d\theta)^n.
\end{split}
\end{equation*}
The coefficient of  $\log(-\epsilon)$ in this expression is 
$$
-2\int_M(n+1)v^{(0)}B\ \theta\wedge (d\theta)^n,
$$
which by (\ref{eq:volpull5}) is
$$
-\frac{(n+1)(-1)^n}{2^nn!}\int_M B \ \theta\wedge (d\theta)^n.
$$
\end{proof}

The definition of CR $Q$-curvature (\cite{fh}) uses the conformal 
structure of Fefferman~\cite{fef76} mentioned in the Introduction. 
As all considerations are local we will assume 
that our strictly pseudoconvex pseudohermitian manifold $(M,\theta)$ is 
 the boundary of a strictly pseudoconvex domain of $\C^{n+1}$. Recall 
from~\cite{fef76} that for a local solution $\rho$ to Fefferman's 
Monge--Amp\'ere problem (\ref{eq:fef}) there is an
associated K\"ahler--Lorentz metric $\gbt$ on $\C^\ast \times X$, namely
$$
\gbt := \sum_{j,k = 0}^{n+1} \frac{\partial^2}{\partial z^j\partial
\overline{z}^k}(|z^0|^2\rho(z))dz^j d\overline{z}^k,
$$
where $(z^0, z) \in \C \times X$. Denote by $\gb$ the restriction of 
$\gbt$ to
$$
\Mb := S^1 \times M = \{(z^0,z) : |z^0| = 1, z\in M \}.
$$
The Lorentz conformal structure $(\Mb,[\gb])$ is known as Fefferman's 
conformal structure. In~\cite{lee86}, Lee showed that corresponding to any 
pseudohermitian structure $\theta$ is an element of the conformal class 
$[\gb]$, call it $\gb(\theta)$, satisfying $\gb(e^{2\Upsilon}\theta) = 
e^{2\Upsilon}\gb(\theta)$.
In this formulation the distinguished representative $\gb$ corresponds to the
pseudohermitian structure $\Im\:\partial \rho|_{TM}$. The $Q$-curvature 
of the conformal structure $(\Mb,[\gb])$ turns out to be $S^1$-invariant and
the \textit{CR $Q$-curvature} $Q^{\textnormal{CR}}$ is defined to be its 
projection down to $M$. 

If we apply Theorem~\ref{thm:fg} (with $\nb = 2(n+1)$) to Fefferman's 
conformal structure, then
$Q^{\textnormal{CR}}$ is a constant multiple of $\Bb|_{S^1\times M}$. To 
relate this to $B|_M$ from Proposition~\ref{prop:lap} we need to compare the 
Poincar\'e  Laplacian $\Db_{\gbp}$ with the K\"ahler Laplacian 
$\Delta_{g_+}$. The Poincar\'e metric $\gbp$ for the 
conformal structure $(\Mb,[\gb])$ lives on $\Xbo := S^1\times X$ and it was 
shown in~\cite{fg85} that
$$
\gbp = g_+ - \left(d\lambda - \frac{i}{2}(\dbar - \partial)u\right)^2,
$$
where $\lambda$ is a coordinate on $S^1$ and $u:= \log(-1/\rho)$. Let
capital indices run from 0 to $2n+2$, with $\jb = j + n + 1$, and  
denote, e.g., $(\partial/\partial z^j) u$ by $u_j$. With this notation, 
$g_+ = \sum u_{j\kb}dz^jd\overline{z}^k$. The matrix of $\gbp$ is
$$
((\gbp)_{JK}) =
\left(\begin{array}{ccc}
\frac{1}{4}u_ju_k & \frac{1}{2}u_{j\kb} - \frac{1}{4}u_j u_{\kb}
& -\frac{i}{2}u_j\\
\frac{1}{2}u_{j\kb} - \frac{1}{4}u_j u_{\kb} & \frac{1}{4}u_{\jb}u_{\kb} & 
\frac{i}{2}u_{\jb}\\
-\frac{i}{2}u_k & \frac{i}{2}u_{\kb} & -1
\end{array}\right).
$$
There is the following matrix factorisation:
$$
((\gbp)_{JK}) = 
\left(\begin{array}{ccc}
1 & 0 & \frac{i}{2}u_j\\
1 & 1 & -\frac{i}{2}u_{\jb}\\
0 & 0 & 1
\end{array}\right)
\left(\begin{array}{ccc}
0 & \frac{1}{2}u_{j\kb} & 0\\
\frac{1}{2}u_{j\kb} & 0 & 0\\
0 & 0 & -1
\end{array}\right)
\left(\begin{array}{ccc}
1 & -1 & 0\\
0 &  1 & 0\\
\frac{i}{2}u_k & -\frac{i}{2}u_{\kb} & 1
\end{array}\right).
$$
From this we can compute the inverse matrix
\begin{equation}
\label{eq:ginv}
((\gbp^{-1})^{JK}) =
\left(\begin{array}{ccc}
0 & 2(u^{-1})^{j\kb} & i(u^{-1})^{j\lb}u_{\lb}\\
2(u^{-1})^{j\kb} &  0 & -i(u^{-1})^{l\jb}u_l\\
i(u^{-1})^{k\lb}u_{\lb} & -i(u^{-1})^{l\kb}u_l & -1 + (u^{-1})^{l\mb}u_lu_{\mb}
\end{array}\right),
\end{equation}
where $(u^{-1})^{j\kb}$ denotes the components of the inverse matrix of 
$(u_{j\kb})$. If $\Deeb$ denotes the Levi-Civita connection of $\gbp$, then 
the Poincar\'e Laplacian is given by
$$
\Db_{\gbp} = -(\gbp^{-1})^{JK}\Deeb_{\partial_{J}}\Deeb_{\partial_{K}}.
$$
On the other hand, the K\"ahler Laplacian is
$$
\Delta_{g_+} = -2(u^{-1})^{j\kb}\partial_j\partial_{\kb}.
$$
After computing Christoffel symbols using the formulae above for $\gbp$ and 
its inverse, it turns out that
\begin{equation}
\label{eq:laplap}
\Db_{\gbp} \equiv 2\Delta_{g_+} \mod \partial_\lambda.
\end{equation}

Using (\ref{eq:ginv}), the definition of $\pb$ in Theorem~\ref{thm:fg} and 
Lee's correspondence mentioned above,
it follows that if $\phi$ is considered as a defining function for 
$S^1\times M$ by extending trivially in the $S^1$ direction, 
then $\pb = -\sqrt{-\phi}$. 
Furthermore, if the function $U$ from Proposition~\ref{prop:lap} is 
considered as a function on $\Xbo$ by extending trivially in the $S^1$
direction, then 
(\ref{eq:laplap}) and Proposition~\ref{prop:lap} imply that 
$$
\frac{1}{2}\Db_{\gbp} U = \nb + O(\pb^{\nb +1}\log(-\pb)). 
$$
By the uniqueness assertion of Theorem~\ref{thm:fg} we conclude that
$$
U = 2\Ubo + O(\pb^{\nb}),
$$
or
\begin{equation*}
\begin{split}
\log(-\phi) + A + B\phi^{n+1}\log(-\phi) &= 2(\log(-\pb) + \Ab 
+ \Bb\pb^{\nb}\log(-\pb)) + O(\pb^{\nb})\\
&\hspace{-8mm}= \log(-\phi) + 2\Ab + \Bb\phi^{n+1}\log(-\phi) 
+ O(\phi^{n+1}).
\end{split}
\end{equation*}
Hence the trivial extension of $B|_M$ to $S^1\times M$ is precisely 
$\Bb|_{S^1\times M}$. We have thus proved (\ref{eq:CRQ}).

\bibliographystyle{amsalpha}
\bibliography{biblio}

\ifx\undefined\bysame
\newcommand{\bysame}{\leavevmode\hbox to3em{\hrulefill}\,}
\fi
\begin{thebibliography}{Web78}

\bibitem[And01]{an}
M.~T. Anderson, {\em ${L}^2$ curvature and volume renormalization of {AHE}
  metrics on 4-manifolds}, Math. Res. Lett. {\bf 8} (2001), 171--188.

\bibitem[BE88]{be88}
D.~Burns and C.~L. Epstein, {\em A global invariant for three dimensional
  {CR}-manifolds}, Invent. Math. {\bf 92} (1988), 333--348.

\bibitem[BE90]{be}
D.~Burns and C.~L. Epstein, {\em Characteristic numbers of bounded domains},
  Acta Math. {\bf 164} (1990), 29--71.

\bibitem[BH05]{bh}
O.~Biquard and M.~Herzlich, {\em A {B}urns-{E}pstein invariant for {ACHE}
  4-manifolds}, Duke Math. J. {\bf 126} (2005), 53--100.

\bibitem[Che45]{chern}
S.~S. Chern, {\em On the curvatura integra in a {R}iemannian manifold}, Ann. of
  Math. {\bf 46} (1945), 674--684.

\bibitem[CY80]{cy}
S.~Y. Cheng and S.~T. Yau, {\em On the existence of a complete {K\"a}hler
  metric on noncompact complex manifolds and the regularity of {F}efferman's
  equation}, Comm. Pure Appl. Math. {\bf 33} (1980), 507--544.

\bibitem[Fef76]{fef76}
C.~Fefferman, {\em Monge--{A}mp{\'e}re equations, the {B}ergman kernel, and the
  geometry of pseudoconvex domains}, Ann. of Math. {\bf 103} (1976), 395--416,
  correction: \textbf{104} (1976), 393--394.

\bibitem[FG85]{fg85}
C.~Fefferman and C.~R. Graham, {\em Conformal invariants}, \'Elie Cartan et les
  math\'ematiques d'aujord'hui, Ast\'erisque, Numero Hors Serie, 1985,
  pp.~95--116.

\bibitem[FG02]{fg}
C.~Fefferman and C.~R. Graham, {\em ${Q}$-curvature and {P}oincar\'e metrics},
  Math. Res. Lett. {\bf 9} (2002), 139--151.

\bibitem[FH03]{fh}
C.~Fefferman and K.~Hirachi, {\em Ambient metric construction of
  ${Q}$-curvature in conformal and {CR} geometries}, Math. Res. Lett. {\bf 10}
  (2003), 819--831.

\bibitem[GL88]{gl}
C.~R. Graham and J.~M. Lee, {\em Smooth solutions of degenerate {L}aplacians on
  strictly pseudoconvex domains}, Duke Math. J. {\bf 57} (1988), 697--720.

\bibitem[Gra00]{graham}
C.~R. Graham, {\em Volume and area renormalizations for conformally compact
  {E}instein metrics}, Rend. Circ. Mat. Palermo (2) Suppl. {\bf 63} (2000),
  31--42.

\bibitem[Her]{herzlich}
M.~Herzlich, {\em A remark on renormalized volume and {E}uler characteristic
  for {ACHE} 4-manifolds}, preprint, {\tt arXiv:math.DG/0305134}.

\bibitem[Hir92]{hirachi}
K.~Hirachi, {\em Scalar pseudo-hermitian invariants and the {S}zeg{\"o} kernel
  on three-dimensional {CR} manifolds}, Complex Geometry, Lecture Notes in Pure
  and Appl. Math., vol. 143, Dekker, 1992, pp.~67--76.

\bibitem[HS98]{hs}
M.~Henningson and K.~Skenderis, {\em The holographic {W}eyl anomaly}, J. High
  Energy Phys. (1998), no.~7, paper 23.

\bibitem[Lee86]{lee86}
J.~M. Lee, {\em The {F}efferman metric and pseudohermitian invariants}, Trans.
  Amer. Math. Soc. {\bf 296} (1986), 411--429.

\bibitem[Lee88]{lee88}
J.~M. Lee, {\em Pseudo-{E}instein structures on {CR} manifolds}, Amer. J. Math
  {\bf 110} (1988), 157--178.

\bibitem[Tan75]{tanaka}
N.~Tanaka, {\em A differential study on strongly pseudo-convex manifolds},
  Lectures in Mathematics, Department of Mathematics, Kyoto University, no.~9,
  Kinokuniya Book-Store Co., Ltd., 1975.

\bibitem[Web78]{webster}
S.~M. Webster, {\em Pseudohermitian structures on a real hypersurface}, J.
  Differential Geom. {\bf 13} (1978), 25--41.

\end{thebibliography}
 
\end{document}